\documentclass[12pt]{article}
\usepackage{amsmath,amsfonts,amssymb,amsthm,graphicx} 

\numberwithin{equation}{section}

\newcommand{\I}{{\bf 1}}
\newtheorem{proposition}{Proposition}[section]
\newtheorem{theorem}[proposition]{Theorem}
\newtheorem{corollary}[proposition]{Corollary}
\newtheorem{lemma}[proposition]{Lemma}
\newtheorem{remark}[proposition]{Remark}

\newtheorem{example}[proposition]{Example}

\newcommand{\nc}{\newcommand}
\nc{\R}{{\mathbb R}}
\nc{\N}{{\mathbb N}}
\nc{\Z}{{\mathbb Z}}

\nc{\BP}{\mathbb{P}}
\nc{\BE}{\mathbb{E}}
\nc{\BQ}{\mathbb{Q}}
\nc{\bN}{{\mathbf N}}
\nc{\BX}{{\mathbb X}}
\nc{\BY}{{\mathbb Y}}
\nc{\bM}{{\mathbf M}}
\nc{\bT}{{\mathbf T}}
\nc{\bi}{{\mathbf 1}}
\nc{\cF}{{\mathcal F}}
\nc{\cK}{{\mathcal K}}
\nc{\cP}{{\mathcal P}}
\nc{\cT}{{\mathcal T}}
\nc{\cS}{{\mathcal S}}
\nc{\cC}{{\mathcal C}}
\DeclareMathOperator{\BV}{{\mathbb Var}}
\DeclareMathOperator{\CV}{{\mathbb Cov}}

\begin{document}

\author{G\"unter Last\footnote{
Institut f\"ur Stochastik, Karlsruher Institut f\"ur Technologie,
76128 Karlsruhe, Germany. 
Email: guenter.last@kit.edu}
\ and \ Eva Ochsenreither\footnote{
Institut f\"ur Stochastik, Karlsruher Institut f\"ur Technologie,
76128 Karlsruhe, Germany. 
Email: e.ochsenreither@kit.edu}  
}

\title{Percolation on stationary tessellations:\\ 
models, mean values and second order structure}
\date{\today}
\maketitle
\begin{abstract}
\noindent
We consider a stationary face-to-face tessellation $X$ of $\R^d$
and introduce several percolation models by colouring
some of the faces black in a consistent way. 
Our main model is cell percolation, where
cells are declared black with probability $p$
and white otherwise. We are interested in geometric
properties of the union $Z$ of black faces.
Under natural integrability assumptions we
first express asymptotic mean-values of intrinsic volumes 
in terms of Palm expectations associated with the faces. 
In the second part of the paper
we study asymptotic covariances of intrinsic volumes of $Z\cap W$,
where the observation window $W$ is assumed to be a polytope. 
Here we need to assume the existence of suitable asymptotic covariances
of the face processes of $X$. We check these assumptions in the
important special case of a Poisson Voronoi tessellation.
In the case of cell percolation on a normal tessellation, 
especially in the plane, our formulae simplify considerably.
\end{abstract}


\section{Introduction}\label{secintro}

Let $X$ be a {\em face-to-face tessellation} of $\R^d$,
that is a random collection of convex and bounded polytopes
(called {\em cells}) covering the whole space and such that for any
different $C,C'\in X$ the intersection $C\cap C'$ is either
empty, or a face of both $C$ and $C'$.
We assume that any bounded subset of $\R^d$ is intersected
by only finitely many cells. 
We interpret $X$ as a {\em point process} on the space of polytopes 
and assume that $X$ is {\em stationary}, meaning that the distribution
of $X$ coincides with that of $\{C+x:C\in X\}$
for all $x\in\R^d$. 
Let, for $k\in\{0,\dots,d\}$, $X_k$ denote the point process
of $k$-dimensional faces of cells in $X$.
We assume throughout that the intensity measure
of $X_k$ is {\em locally finite}. For more details on stationary
tessellations we refer to \cite[Chapter 10]{SW}
and the next section.

For $p\in[0,1]$ and $n\in\{0,\dots,d\}$
we define {\em $n$-percolation} on $X$ as follows.
Given $X$, we colour the polytopes in $X_n$ independently {\em black}
with probability $p$. All other polytopes in $X_n$ are
{\em white}. If $n\le d-1$ and $k\in\{n+1,\dots,d\}$, then we colour $F\in X_{k}$
black if all its $(k-1)$-faces are black.
We are interested in the union $Z$ of all
black faces of $X$. This is a stationary random closed set,
see \cite[Chapter 2]{SW}. In the case $n=d$ we refer to this
as \emph{cell percolation} and for $n=0$ as \emph{vertex percolation}. 
In the planar case $d=2$ we refer to 1-percolation 
as \emph{edge percolation}. In the general case we also speak
of \emph{face percolation}.

\begin{figure}[h!]\label{Abb}
\centering
\includegraphics{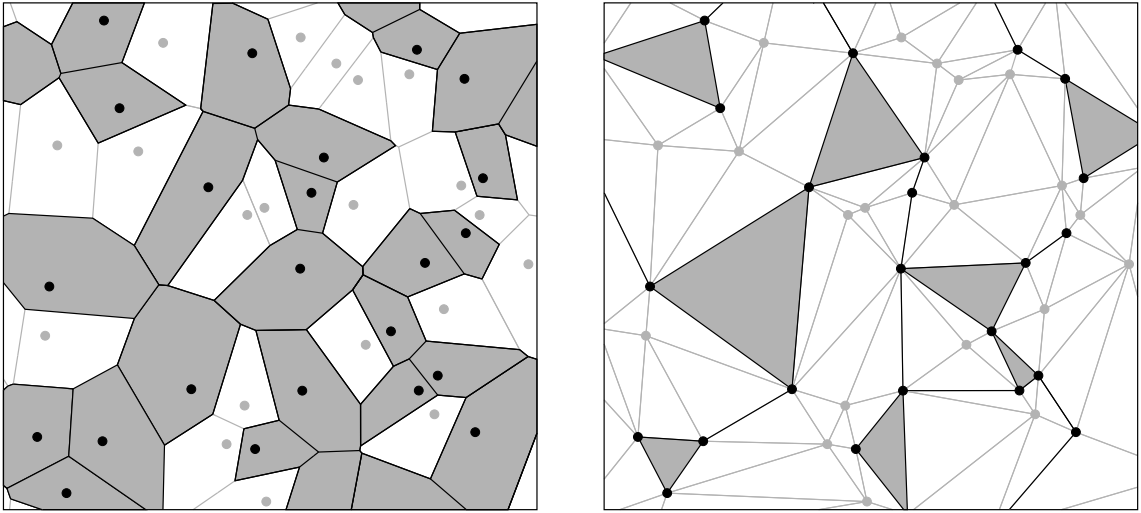}
\caption{Cell percolation on a Poisson Voronoi tessellation 
and vertex percolation on a Poisson Delaunay tessellation.}
\end{figure}

Cell percolation on a Poisson Voronoi tessellation, see Figure \ref{Abb}, 
was studied
in \cite{BR06a}, where it is shown   
that the {\em critical probability} of this model of 
{\em continuum percolation} is $1/2$. The present paper
was motivated by \cite{neher08}, introducing vertex, edge and cell
percolation on several planar lattices. The authors of \cite{neher08} notice that
in many models the only non-trivial zero of the mean Euler
characteristic is a remarkable accurate approximation
of the critical probability. 

Our first aim in this paper is to establish $n$-percolation 
on $X$ as an interesting model of stochastic geometry and continuum percolation. 
Our main aim is to study first and second order geometric properties 
of the black phase $Z$.
Section \ref{secnotation} collects some preliminaries on
stationary tessellations and Palm probabilities
and gives the definition of face percolation.
Asymptotic mean values of intrinsic volumes of $Z\cap W$
are studied in  Sections \ref{secmean} and \ref{meansplanar},
where we assume that the {\em observation window} $W$ is a convex polytope.
Asymptotic covariances of intrinsic volumes are treated
in Sections \ref{secondnormal} assuming that
the asymptotic covariances of intrinsic volumes of face processes
exist. Theorem \ref{t5.1} shows
that these covariances are polynomials in the colouring probability $p$,
where the coefficients are determined
by the global fluctuation of the intrinsic volumes within
the face processes as well as by the local geometry of $X$.
The important special case of cell percolation on a
planar and {\em normal} tessellation is discussed
in Section \ref{secvEuler}.  
In Section \ref{poissonvoronoi} we check that
a Poisson Voronoi tessellation satisfies all assumptions
required for our general results.
Moreover, all asymptotic covariances are then given by
fairly explicit integral formulae.  
For cell percolation in the planar case, the asymptotic variance of the Euler
characteristic is determined by the intensity and the second
moment of the number of vertices of a typical cell and has
a global maximum at the critical threshold $p=1/2$, see Corollary \ref{varianceEuler}.
The Appendix contains some integrability properties of 
a Poisson Voronoi tessellation.

\section{Notation and preliminaries}\label{secnotation}

\subsection{Palm calculus}\label{secPalm}

It is convenient to follow \cite{La10,LaTho09} by assuming the
basic sample space $(\Omega,\mathcal{F})$ to 
be equipped with a {\em measurable flow}
$\theta_x:\Omega \to \Omega$, $x\in \R^d$, that is
$(\omega,x)\mapsto \theta_x\omega$ is measurable,
$\theta_{x+y}=\theta_x\circ\theta_y$ for all $x,y\in\R^d$ and
$\theta_0$ is the identity on $\Omega$. We further assume
that $\BP$ is {\em stationary}, i.e.\ $\BP\circ\theta^{-1}_x=\BP$,
$x\in\R^d$. 
A \emph{random measure} $\mu$ on $\R^d$ is a kernel from $\Omega$ to $\R^d$, such that
$\mu(\omega,\cdot)$ is locally finite for all $\omega\in\Omega$.
If $\mu(\omega,B)$ is integer-valued for all bounded Borel sets $B\subset\R^d$,
then $\mu$ is a {\em point process}, that is called {\em simple}
if $\mu(\{x\})\le 1$ for all $x\in\R^d$. In the latter case 
we identify $\mu$ with its support $\{x\in\R^d:\mu(\{x\})>0\}$.
A random measure $\mu$ is called {\em invariant} if
\begin{align}\label{inv}
\mu(\theta_x\omega,B-x)=\mu(\omega,B),\quad x\in\R^d,\,\omega\in\Omega,
\end{align}
for any Borel set $B\subset\R^d$. It then follows that
$\mu$ is {\em stationary}, that is the distribution
of $\mu(\cdot +x)$ is independent of $x\in\R^d$.
If $\mu$ is invariant, then
$\gamma_\mu:=\BE \mu([0,1]^d)$ is called
{\em intensity} of $\mu$. If 
$0<\gamma_\mu<\infty$,
then the {\em Palm probability measure} $\BP^0_\mu$ of $\mu$ is
defined by 
\begin{align} \label{Palm}
\BP^0_\mu(A):=\gamma_\mu^{-1}\iint \I_A(\theta_x\omega)\I\{x\in[0,1]^d\}\,
\mu(\omega,\mathrm{d}x)\,\BP(\mathrm{d}\omega), \quad A\in\mathcal{F}.
\end{align}
It satisfies the {\em refined Campbell theorem}
\begin{align} \label{refCampbell}
\BE \int f(\theta_x,x)\,\mu(\mathrm{d}x)=
\gamma_\mu\BE^0_\mu\int f(\theta_0,x)\,\mathrm{d}x
\end{align}
for all measurable $f:\Omega\times\R^d\to [0,\infty)$, where
$\BE^0_\mu$ denotes expectation with respect to
$\BP^0_\mu$.

For ease of reference we now state {\em Neveu's exchange formula}.
It will be frequently used in this paper. 
This formula
also goes under the name {\em mass-transport principle},
see \cite{La10,LaTho09} for a brief discussion.

\begin{proposition}\label{mtp}
Let $\mu$ and $\mu'$ be invariant random measures on $\R^d$
with positive and finite intensities
and let $h:\Omega\times\R^d\to[0,\infty)$
be measurable. Then
\begin{align}\label{neveu3}
\gamma_\mu\BE^0_{\mu}\bigg[\int h(\theta_x,-x)\,\mu'(\mathrm{d}x)\bigg]=
\gamma_{\mu'}\BE^0_{\mu'}\bigg[\int h(\theta_0,x)\,\mu(\mathrm{d}x)\bigg].
\end{align}
\end{proposition}

\subsection{Coloured tessellations and face percolation}

We start with introducing some basic terminology for tessellations
and refer to \cite{SW} for further detail.
We let $\mathcal{K}^d$ denote the space of convex bodies
(convex and compact subsets of $\R^d$) and equip it
with the Borel $\sigma$-field
associated with the {\em Hausdorff distance}.
A {\em polytope} is a finite
intersection of half-spaces which is bounded and non-empty.
The system $\cP^d$ of all such polytopes is a measurable
subset of $\mathcal{K}^d$. 
A {\em tessellation} (of $\R^d$) is a countable system $\varphi$ of
polytopes (cells) covering the whole space such that any two
different elements of $\varphi$ have disjoint interior and
any bounded subset of $\R^d$ is intersected
by only finitely many cells. 
Let $k\in\{0,\dots,d-1\}$. A {\em $k$-face}
of $C\in\cP^d$ is a $k$-dimensional intersection of $C$
with a {\em supporting hyperplane} of $C$. We let $\cF_k(C)$
denote the system of all $k$-faces of $C$. It is convenient
to define $\cF_d(C):=\{C\}$.
A tessellation $\varphi$ is {\em face-to-face}
if for $C,C'\in \varphi$ the intersection $C\cap C'$ is either
empty, or a face of both $C$ and $C'$.
Let $\bT$ denote the set of all face-to-face tessellations.
We define the system of $k$-faces of $\varphi\in\bT$ by
\begin{align}\label{2.1}
\cF_k(\varphi):=\bigcup_{C\in\varphi}\cF_k(C)
\end{align}
and the system of faces of $\varphi$ by
\begin{align}\label{2.1a}
\cF(\varphi):=\bigcup_{k=0}^d\cF_k(\varphi).
\end{align}
Note that $\cF_d(\varphi)=\varphi$.

In this paper we define a
{\em coloured tessellation} as a tuple $\psi=(\varphi,\varphi_0,\dots,\varphi_d)$, 
where $\varphi\in\bT$ and $\varphi_k\subset \cF_k(\varphi)$ such that 
$\cF_{k-1}(F)\subset\varphi_{k-1}$ whenever $k\ge 1$ and
$F\in\varphi_k$.
Any face in $\cup^d_{k=0}\varphi_k$ is called
{\em black}, while the other faces of $\varphi$ are called {\em white}. 
If $F \in \cF(\varphi)$ is black, then by definition all its faces are black as well.
We write $X(\psi):=\varphi$ and $X_k^1(\psi):=\varphi_k$.

Let $\bT_c$ denote the space of all coloured tessellations.
We identify discrete sets with the associated counting measures.
In particular we write for measurable $H\subset\cP^d$ and $\psi \in \bT_c$
\begin{align}\label{2.4a}
X(\psi,H)&:=|\{C\in X(\psi):C\in H\}|,\\
\label{2.4b}
X_k^1(\psi,H)&:=|\{F\in X_k^1(\psi):F\in H\}|,
\end{align}
where $|A|$ denotes the cardinality of a set $A$.
Let $\cT_c$ denote the smallest $\sigma$-field on
$\bT_c$ so that $\psi\mapsto (X(\psi,H),X_0^1(\psi,H_0),\dots,X_d^1(\psi,H_d))$
is measurable for all measurable $H,H_0,\dots,H_d\subset\cP^d$. 
The $\sigma$-field  $\cT$ on $\bT$ is defined similarly.

A {\em random coloured tessellation} $\Psi$ is a measurable mapping
from the probability space $(\Omega,\mathcal{F},\BP)$
to $(\bT_c,\cT_c)$. 
In particular, $X(\Psi)$ and $X_0^1(\Psi),\dots,X_d^1(\Psi)$
are then point processes on $\cP^d$, see also \eqref{2.4a} and \eqref{2.4b}. 
The same is true for  $\mathcal{F}_0(X(\Psi)),\dots,\mathcal{F}_d(X(\Psi))$.
We are interested in the union
\begin{align}\label{Z}
  Z:=\bigcup^d_{k=0}\bigcup_{F\in X_k^1(\Psi)} F
\end{align}
of all black faces. It can be shown, that $Z$ is a {\em random closed set},
see \cite{SW} for a definition of this concept.
We shall always assume that $\Psi$ is stationary, that is,
\begin{align}\label{2.3}
\Psi+x\overset{d}{=}\Psi,\quad x\in\R^d,
\end{align} 
where, for $\psi=(\varphi,\varphi_0,\dots,\varphi_d)\in\bT_c$, 
$\psi+x:=(\varphi+x,\varphi_0+x,\dots,\varphi_d+x)$,
$H+x:=\{F+x:F\in H\}$ for $H\subset\cK^d$, and $A+x:=\{y+x:y\in A\}$
for $A\subset \R^d$. In that case $Z$ is stationary as well, that is,
\begin{align}\label{2.4}
Z+x\overset{d}{=}Z,\quad x\in\R^d.
\end{align}

We will be mainly concerned with what we call $n$-percolation
(or face percolation) on a stationary tessellation. To introduce this concept
we assume given a random face-to-face tessellation $X$,
that is a random element of the space $\bT$. We assume that
$X$ is stationary, i.e.\ that the distribution of $X+x$
does not depend on $x\in\R^d$. A coloured tessellation
$\Psi$ is an \emph{$n$-percolation} on $X$ with (percolation) parameter $p$ if $X(\Psi)=X$,
the point process $X_n^1(\Psi)$ is a
{\em $p$-thinning} of $\cF_n(X)$ (see \cite{Kallenberg} for a definition
of a thinning) and if $X_0^1(\Psi),\dots,X_d^1(\Psi)$
are given in the following way. For $k<n$ the system $X_k^1(\Psi)$ is the union
of all $\cF_k(F)$ for $F\in X_n^1(\Psi)$. For $k>n$ the
system $X_k^1(\Psi)$ is defined recursively. A polytope $F\in X_k$
belongs to $X_k^1(\Psi)$ if and only if $\cF_{k-1}(F)\subset X_{k-1}^1(\Psi)$.
In the case $n=d$ we speak of {\em cell percolation} and
in the case $n=0$ of {\em vertex percolation}.

We now fix a coloured tessellation $\Psi$  such that
\begin{align}\label{adapt}
\Psi(\theta_x\omega)=\Psi(\omega)-x,\quad \omega\in \Omega,\,x\in \R^d.
\end{align}
Then $\Psi$ is stationary in the sense of \eqref{2.3}.
Throughout we will use the
following short-hand notation
for the systems of all (respectively all black) $k$-faces:
\begin{align}\label{3.1}
X_k:=\mathcal{F}_k(X),\quad X_k^1:=X_k^1(\Psi),\quad k\in\{0,\dots,d\}.
\end{align}
The invariance assumption \eqref{adapt} implies
\begin{align}\label{3.11}
(X_k(\theta_x\omega),X_k^1(\theta_x\omega))
=(X_k(\omega)-x,X_k^1(\omega)-x),\quad \omega\in \Omega,\,x\in \R^d.
\end{align}
For $k\in\{0,\dots,d\}$ let
\begin{align}\label{etak}
\eta^{(k)}:=\{s(F):F\in X_k\}
\end{align}
be the point process of
{\em Steiner points} of the faces in $X_k=\cF_k(X)$,
see \cite{SW} for the definition of the Steiner point $s(K)$
of a non-empty $K\in\cK^d$. Since $s(K+x)=s(K)+x$ for all
$x\in\R^d$, \eqref{3.11} implies that $\eta^{(k)}$ is invariant.
By assumption on $X$, $\eta^{(k)}$ contains infinitely many points
so that the intensity
\begin{align}\label{ietaj}
\gamma_k:=\gamma_{\eta^{(k)}}=\BE\eta^{(k)}([0,1]^d)
\end{align}
is positive. We assume $\gamma_k<\infty$, so that
the Palm probability measure
$\BP^0_k:=\BP^0_{\eta^{(k)}}$ is well-defined.
The expectation
with respect to $\BP^0_k$ is denoted by $\BE^0_k$.
Note that under $\BP^0_k$ the origin is almost surely in the relative interior
of some $k$-dimensional face.

Let $\psi=(\varphi,\varphi_0,\dots,\varphi_d)$ be a coloured tessellation
and let $x\in\R^d$. Since $\varphi$ is face-to-face, there is
unique $F\in\mathcal{F}(\varphi)$ such that
$x$ is in the relative interior of $F$. We then write
$F(\psi,x)\equiv F(\varphi,x) = F$. 
To treat the local neighbourhood of $x\in\R^d$ we introduce,
for $l\in\{0,\dots,d\}$, the set $\cS_l(\psi,x)\equiv \cS_l(\varphi,x)$
as follows. Let $k$ be the dimension of  $F(\psi,x)$.
If $l\geq k$ (resp.\ $l<k$) then we let $\cS_l(\psi,x)$ be the set of
all faces $G\in\mathcal{F}_l(\varphi)$ such that
$F(\psi,x)\subset G$ (resp.\ $G\subset F(\psi,x)$).
It is convenient to abbreviate
\begin{align*}
(F(x),\cS_l(x)):=(F(\Psi,x),\cS_l(\Psi,x)),\quad x\in\R^d.
\end{align*}
Since $F(\Psi,x)=F(\Psi-x,0)+x$ we obtain from
\eqref{adapt} that
\begin{align}\label{3.6} 
\BP^0_k(F(0)\in\cdot)=\gamma^{-1}_k\BE \int\I\{x\in [0,1]^d,F(x)-x\in\cdot\}\,\mathrm{d}x
\end{align}
is the distribution of a {\em typical $k$-face}.
The next result is a version of Theorem 10.1.1 in \cite{SW}.
The proof can easily be given with Neveu's exchange formula, see also \cite{BaumLa07}.

\begin{proposition}\label{t2.1} Let $k,l\in\{0,\dots,d\}$
and $g:\cP^d\times\cP^d\to[0,\infty)$ be a
measurable function. Then
\begin{align}\label{2.17}
\gamma_k\BE^0_k\sum_{G\in\cS_l(0)}g(F(0),G-s(G))
=\gamma_l\BE^0_l\sum_{F\in\cS_k(0)}g(F-s(F),F(0)).
\end{align}
\end{proposition}

In particular, Proposition \ref{t2.1} implies that
\begin{align}\label{2.18}
\gamma_k n_{k,l}=\gamma_l n_{l,k},\quad k,l\in\{0,\dots,d\},
\end{align}
where 
\begin{align}\label{2.19}
n_{k,l}:=\BE^0_k|\cS_l(0)|.
\end{align}
We refer to Section 10.1 of \cite{SW} for further information
on such {\em face star} relationships.

\section{Mean value analysis}\label{secmean}

Let $X$ be a stationary face-to-face tessellation,
that is, a random element of $\bT$. Let, for $k\in\{0,\dots,d\}$, 
$X_k=\cF_k(X)$ denote the point process of $k$-faces of $X$.
We assume that 
\begin{align}\label{locfinite}
\sum^d_{k=0}\BE\sum_{F\in X_k}\I\{F\cap K\ne\emptyset\}<\infty,\quad K\in\cK^d,
\end{align}
an assumption, quite common in stochastic geometry \cite{SW}.
It is easy to see that \eqref{locfinite} implies $\gamma_k:=\BE \eta^{(k)}([0,1]^d)<\infty$,
where the point process $\eta^{(k)}$ is defined by \eqref{etak}.
The refined Campbell theorem \eqref{refCampbell} allows to
rewrite \eqref{locfinite} as
\begin{align}\label{locfinite3}
\sum^d_{k=0}\BE^0_k V_d(F(0)+K)<\infty,\quad K\in\mathcal{K}^d,
\end{align}
where we have used that $(A+x)\cap B\ne\emptyset$
for $A,B\subset\R^d$ and $x\in\R^d$ iff $x\in B-A:=\{y-z:y\in A,z\in B\}$. 
Recall that $F(x) \in \cF(X)$ is the unique face that contains $x \in \R^d$ 
in its relative interior.
Often we have to assume that
\begin{align}\label{locfinite2}
\sum^d_{i,k=0}\BE^0_k V_i(F(0))^2<\infty.
\end{align}
Note that \eqref{locfinite3} is a consequence of \eqref{locfinite2}, 
the Steiner formula and the Cauchy-Schwarz inequality.

For $n\in\{0,\dots,d\}$ we consider $n$-percolation $\Psi$ on $X$.
It is no restriction of generality to assume that \eqref{adapt} holds.
Let the stationary random closed set $Z$ be given by \eqref{Z}.
The density of the $i$-th intrinsic volume  of $Z$
is defined by the limit
\begin{align}\label{621}
\delta_i(p):=\lim_{t\to \infty}V_d(W_t)^{-1}\BE V_i(Z\cap W_t),
\end{align}
where $W_t:=t^{1/d}W$ and $W\in\mathcal{P}^d$ 
is assumed to have volume one and to contain the origin in its interior.
We shall show below that this limit exists and does not depend on $W$.  
Our first aim in this paper is to derive a
polynomial formula for these densities.
It should not come as surprise, that this formula is based
on the joint distribution of $(V_i(F(0)),|\mathcal{S}_n(0)|)$
under the measures $\BP^0_k$.

\begin{theorem}\label{t1a} Consider $n$-percolation on $X$ and let $i\in\{0,\dots,d\}$. 
Assume \eqref{locfinite2}. Then the limit \eqref{621} exists and is given by
\begin{align}\label{densities}
\delta_i(p)=&\sum^{n-1}_{k=i}(-1)^{i+k}\gamma_k 
\BE^0_k\big[(1-(1-p)^{|\mathcal{S}_n(0)|})V_i(F(0))\big]\\
& + \sum_{k=n}^d (-1)^{i+k}\gamma_k 
\BE^0_k\big[p^{|\mathcal{S}_n(0)|}V_i(F(0))\big],\notag \quad p\in[0,1].
\end{align}
In particular, for cell percolation we have
\begin{align}\label{celldensities}
\delta_i(p)=
\begin{cases}
\sum_{k=i}^d (-1)^{i+k+1} \gamma_k \BE_k^0[(1-p)^{|\cS_d(0)|} V_i(F(0))], \quad i<d,\\
p, \quad i=d.
\end{cases}
\end{align}
\end{theorem}


We prepare the proof of  Theorem \ref{t1a} with some geometric
preliminaries. The intrinsic volumes can be defined for convex bodies 
by the Steiner formula. By additivity they can then be extended to finite 
unions of convex bodies, see e.g.\ \cite{SW}.
Groemer \cite{Groemer72} defines
the intrinsic volumes for a much wider class 
of \emph{approximable sets} 
containing the relative interior of convex bodies and the intersection 
of a relative open polytope with the boundary of a convex body such that 
they are still additive and rigid motion invariant. In particular,
\begin{align}\label{relint}
V_i(\mathrm{relint}(K)) = (-1)^{i+\dim(K)} V_i(K), \quad K \in \mathcal{K}^d,
\end{align}
where $\mathrm{relint}(B)$ denotes the relative interior of a set $B$. 
Let $\mathrm{int}(B)$ and $\partial B$ 
denote the interior and the boundary of a convex body. We can 
write $Z \cap W_t$ as a disjoint union
\begin{align*}
Z \cap W_t &= (Z \cap \mathrm{int}(W_t)) \cup (Z \cap \partial W_t)\\
&= \bigcup_{k=0}^d \bigcup_{F \in X_k^1} (\mathrm{relint}(F) \cap \mathrm{int}(W_t)) 
\cup \bigcup_{k=0}^d \bigcup_{F \in X_k^1} (\mathrm{relint}(F) \cap \partial W_t).
\end{align*}
Since the tessellation $X$ is stationary, the intersection of a $k$-face $F$ 
with $W_t$ is almost surely empty if 
$\mathrm{relint}(F) \cap \mathrm{int}(W_t) = \emptyset$. 
Thus, $\mathrm{relint}(F) \cap \mathrm{int}(W_t) =\mathrm{relint}(F \cap W_t)$ 
a.s. It follows that both $Z \cap \mathrm{int}(W_t)$ and $Z \cap \partial W_t$
are approximable. The additivity of the intrinsic volumes and \eqref{relint} 
yield almost surely
\begin{align}\label{FormelVi}
V_i(Z \cap \mathrm{int}(W_t)) 
&= \sum_{k=0}^d \sum_{F \in X_k^1}V_i(\mathrm{relint}(F \cap W_t)) 
=\sum_{k=0}^d (-1)^{i+k} \sum_{F \in X_k^1} V_i(F \cap W_t)
\end{align}
because $\dim(F \cap W_t)=k$ almost surely for $F \in X_k$ if the intersection 
is non-empty. Since the observation window $W$ is a polytope, 
we can partition $\partial W$ in the relative interior of the lower-dimensional 
faces of $W$ and get
\begin{align*}
V_i(Z \cap \partial W_t) &= \sum_{k=0}^d \sum_{l=0}^{d-1} \sum_{U \in \cF_l(W)} 
\sum_{F \in X_k^1} V_i(\mathrm{relint}(F) \cap \mathrm{relint}(U_t))\\
&= \sum_{k=0}^d \sum_{l=0}^{d-1} \sum_{U \in \cF_l(W)} \sum_{F \in X_k^1} 
V_i(\mathrm{relint}(F \cap U_t)),
\end{align*}
where $U_t:=t^{1/d}U$ denotes similarly as before the scaled face and the last equation 
holds almost surely because of the stationarity of the tessellation since 
the intersection of $F$ and $U_t$ is almost surely empty if 
$\mathrm{relint}(F) \cap \mathrm{relint}(U_t)=\emptyset$. Using \eqref{relint}, 
it follows that
\begin{align}\label{FormelRand}
V_i(Z \cap \partial W_t) = \sum_{k=0}^d \sum_{l=0}^{d-1} \sum_{U \in \cF_l(W)} 
\sum_{F \in X_k^1} (-1)^{i+\dim(F \cap U_t)} V_i(F \cap U_t)
\end{align}
almost surely. 

\bigskip

{\em Proof of Theorem \ref{t1a}.}
First, we will show that 
\begin{align}\label{VorEWRand}
\lim_{t \to \infty} t^{-1} \BE[V_i(Z \cap \partial W_t)]=0.
\end{align}
Because of the representation \eqref{FormelRand} it is enough to show that
\begin{align*}
\lim_{t \to \infty} t^{-1} \BE \bigg[ \sum_{F \in X_k^1} (-1)^{i+\dim(F \cap U_t)} 
V_i(F \cap U_t) \bigg] =0
\end{align*}
for $k \in \{0,\ldots,d\}$ and $U \in \cF(W)$ with $\dim(U)<d$. By the 
definition of $n$-percolation we have
\begin{align}\label{EWViRand}
&\BE \bigg[ \sum_{F \in X_k^1} (-1)^{i+\dim(F \cap U_t)} V_i(F \cap U_t) \bigg]\\
&=
\begin{aligned}[t]
\sum_{r=1}^\infty  &((1-(1-p)^{r}) \I\{k<n\} + p^{r} \I\{k \geq n\})\\
&\times \BE \int (-1)^{i+\dim(F(x) \cap U_t)} V_i(F(x) \cap U_t) \I\{|\cS_n(x)|=r\} \, 
\eta^{(k)}(\mathrm{d}x).
\end{aligned} \notag
\end{align}
Using the monotonicity of the intrinsic volumes, we obtain that
\begin{align*}
\bigg| \BE \bigg[\sum_{F \in X_k^1} &(-1)^{i+\dim(F \cap U_t)} V_i(F \cap U_t) \bigg] \bigg| 
\le \BE \int V_i(F(x) \cap U_t) \, \eta^{(k)}(\mathrm{d}x)\\
&\le \BE \int V_i(F(x)) \I\{F(x) \cap \partial W_t \neq \emptyset\} \, 
\eta^{(k)}(\mathrm{d}x)\\
&=\gamma_k \BE_k^0 \int V_i(F(0))\I\{(F(0)+x)\cap\partial W_t\ne\emptyset\}\,
\mathrm{d}x,
\end{align*}
where we have used the refined Campbell theorem \eqref{refCampbell}
to get the final identity.
We claim that
\begin{align}\label{claim}
\lim_{t \to \infty} \frac{1}{t} \BE \int V_i(F(x)) \I\{F(x) \cap \partial W_t \neq \emptyset\} 
\, \eta^{(k)}(\mathrm{d}x) =0, \quad k \in \{0,\ldots,d\}.
\end{align}
Indeed, we have $\lambda_d(\partial W-t^{-1/d}K)\to \lambda_d(\partial W)=0$
as $t\to\infty$
for any $K\in\mathcal{K}^d$, where $\lambda_d$ denotes the Lebesgue measure on $\R^d$, 
see the proof of Theorem 4.1.3 in \cite{SW}. Moreover, as in the cited proof
we have $\lambda_d(\partial W-t^{-1/d}K)\le c\lambda_d(B^d+K)$
for all $t\ge 1$ and all convex bodies $K$, where $B^d$ is the unit ball
and $c>0$ does not depend on $K$.
Hence \eqref{claim} follows from the Steiner formula, 
the Cauchy-Schwarz inequality, our assumption \eqref{locfinite2} and dominated convergence.
In particular, \eqref{VorEWRand} holds. 


Now we treat the main terms \eqref{FormelVi}. 
The definition of $n$-percolation implies that 
\begin{align}\label{3.13}
\BE V_i(Z \cap \mathrm{int}(W_t))
=
\begin{aligned}[t]
\sum^d_{k=0}\sum^\infty_{r=1}&((1-(1-p)^r) \I\{k<n\} + p^r \I\{k \geq n\})\\
&\times \BE\int V_i(F(x) \cap W_t) \I\{|\cS_n(x)|=r\}\,\eta^{(k)}(\mathrm{d}x).
\end{aligned}
\end{align}
Since, for $k\in\{0,\dots,d\}$ and $x\in\eta^{(k)}$,
\begin{align*}
|V_i(F(x) \cap W_t)- V_i(F(x)) \I\{x \in W_t\}|\le 
V_i(F(x)) \I\{F(x)\cap \partial W_t\ne\emptyset\},
\end{align*}
\eqref{claim} implies that
\begin{align}\label{3.16}
\lim_{t\to\infty} \frac{1}{t} \BE V_i(Z \cap \mathrm{int}(W_t))
= \lim_{t\to\infty}\frac{1}{t} &\sum_{k=0}^d \sum^\infty_{r=1} ((1-(1-p)^r) \I\{k<n\} 
+ p^r \I\{k \ge n\}) \notag\\
& \times \BE\int V_i(F(x)) \, \I\{x\in W_t,|\cS_n(x)|=r\}\,\eta^{(k)}(\mathrm{d}x).
\end{align}
But the refined Campbell theorem \eqref{refCampbell} yields that
\begin{align*}
\frac{1}{t} \BE\int V_i(F(x)) \, \I\{x\in W_t,|\cS_n(x)|=r\}\,\eta^{(k)}(\mathrm{d}x)
=\gamma_k\BE^0_k \big[V_i(F(0)) \I\{|\cS_n(0)|=r\}\big].
\end{align*}
Combining \eqref{3.16} with \eqref{VorEWRand}, yields 
the assertion \eqref{densities}.

For cell percolation, \eqref{densities} equals
\begin{align*}
\delta_i(p) = \sum_{k=i}^d (-1)^{i+k} \, \gamma_k \, \BE_k^0[(1-(1-p)^{|\cS_d(0)|}) V_i(F(0))].
\end{align*}
Applying Theorem 10.1.4 in \cite{SW} gives the second assertion for $i<d$. 
In the case $i=d$ we have
\begin{align*}
\delta_d(p) = (1-(1-p)^1) \, \gamma_d \, \BE_d^0[V_d(F(0))] = p
\end{align*}
since $\gamma_d \, \BE_d^0[V_d(F(0))]=1$.
\qed

\bigskip

The tessellation $X$ is {\em normal} if for $0\le k\le d$
any $k$-face is almost surely contained in $d-k+1$ cells.
In this case we have the following duality relation:

\begin{proposition}\label{BezEW}
Consider cell percolation on a normal tessellation $X$ and assume \eqref{locfinite2}. 
Then we have for $p \in [0,1]$  
\begin{align*}
  \delta_i(p)=(-1)^{d+i+1} \delta_i(1-p), \quad i\in\{0,\dots,d-1\}.
\end{align*}
\end{proposition}
\emph{Proof.} To make the dependence on the colouring probability $p \in [0,1]$
explicit, we write $Z_p$ instead of $Z$. The very definition of face percolation
yield
\begin{align}\label{Z(1-p)}
  \overline{Z_p^\mathrm{c}} \overset{d}{=} Z_{1-p},
\end{align}
where $B^\mathrm{c}$ and $\bar{B}$ denote the complement 
and the closure of a set $B \subset \R^d$, respectively.
Define the set of all white $k$-faces by $X_k^0:= X_k \setminus X_k^1$. 
The additivity of the intrinsic volumes and \eqref{relint} yield almost surely
\begin{align*}
V_i(Z_p^\mathrm{c} \cap \mathrm{int}(W_t)) 
=\sum_{k=0}^d \sum_{F \in X_k^0} V_i(\mathrm{relint}(F) \cap \mathrm{int}(W_t)) 
=\sum_{k=0}^d \sum_{F \in X_k^0} (-1)^{i+k} V_i(F \cap W_t),
\end{align*}
because we have a.s.\ $F \cap W_t = \emptyset$ if 
$\mathrm{relint}(F) \cap \mathrm{int}(W_t) = \emptyset$ and 
$\dim(F \cap W_t)=\dim(F)$ if $F \cap W_t \neq \emptyset$. Since $X$ 
is normal it follows from the inclusion-exclusion principle that
\begin{align*}
V_i(Z_p^\mathrm{c} \cap \mathrm{int}(W_t)) 
=(-1)^{d+i} V_i(\overline{Z_p^\mathrm{c}} \cap W_t).
\end{align*}
Since $V_i(Z_p^\mathrm{c} \cap \mathrm{int}(W_t))+V_i(Z_p \cap \mathrm{int}(W_t))
=V_i(\mathrm{int}(W_t))$ we obtain that
\begin{align}\label{Z(1-p)2}
V_i(\overline{Z_p^\mathrm{c}} \cap W_t) 
=(-1)^{d+i}t^{i/d}V_i(W)+(-1)^{d+i+1} V_i(Z_p \cap W_t),
\end{align}
where we have also used the homogeneity of intrinsic volumes.
Combining this with \eqref{Z(1-p)} and using \eqref{VorEWRand}
(for $Z_{1-p}$) yields the assertion.\qed

\bigskip

Combining Proposition \ref{BezEW} with Theorem \ref{t1a}
we obtain the following result.

\begin{proposition}\label{c1} Consider cell percolation on a normal tessellation $X$. 
Let $i \in \{0,\ldots,d\}$ and assume \eqref{locfinite2}. Then,
\begin{align}\label{meannormal}
\delta_i(p) = \sum_{k=i}^d (-1)^{d-k} p^{d-k+1} \, \gamma_k \BE_k^0[V_i(F(0))].
\end{align}
\end{proposition}

\bigskip

\section{Mean values in the planar case}\label{meansplanar}

In this section we discuss the results of the previous section
in the planar case $d=2$. We assume given a stationary tessellation
$X$ satisfying \eqref{locfinite2}. We start with proving
(the well-known \cite{SW}) equations
\begin{align}\label{vertices}
\gamma_0=\gamma_2\frac{2}{n_{0,1}-2},\qquad
\gamma_1=\gamma_2\frac{n_{0,1}}{n_{0,1}-2},
\end{align}
where we recall from \eqref{2.19} that $n_{0,1}=n_{0,2}$ is the mean degree of a typical
vertex. In particular,
\begin{align}
\gamma_1=\gamma_0+\gamma_2.
\end{align}

\begin{proposition}\label{pint} Assume that $n_{0,1}<\infty$. 
Then the intensities $\gamma_0,\gamma_1$
are given by \eqref{vertices}. Moreover,
\begin{align}\label{meanedges}
n_{2,0}=\frac{2n_{0,1}}{n_{0,1}-2}.
\end{align}
\end{proposition}
{\em Proof.}  Take in \eqref{2.17} $k=0$, $l=2$ and $g(F,G)$ 
as the interior angle of $G$ at $F$ (normalized such
that a full angle equals $1$) if $\dim G=2$ and $F$ is a
vertex of $G$. This yields $\gamma_0=\gamma_2 (n_{2,0}-2)/2$.
Together with $\gamma_2 n_{2,0}=\gamma_0 n_{0,2}$ (see \eqref{2.18}) and
$n_{0,1}=n_{0,2}$ this gives \eqref{meanedges} and the
first equation in \eqref{vertices}. The second equation
can be obtained from \eqref{meanedges} and 
$\gamma_2 n_{2,1}=\gamma_1 n_{1,2}$, that is
$\gamma_2 n_{2,0}=2\gamma_1$.\qed

\begin{remark}\label{r1}\rm If $X$ is normal, then
$\gamma_0=2\gamma_2$, $\gamma_1=3\gamma_2$, and 
$n_{2,0}=6$, see also \cite[Theorem 10.1.6]{SW}.
\end{remark}

Next we state formulae for the asymptotic mean of the Euler characteristic in 
the planar case, using Theorem \ref{t1a} and generalizing 
the results in \cite[Section 2]{neher08}. We write
\begin{align}\label{4.12}
p_{k,n}(m):=\BP^0_k(|\cS_n(0)|=m),\quad k\in\{0,\dots,d\},\, m\in\N.
\end{align}

\begin{corollary}\label{c} Consider cell percolation on
$X$. Then
\begin{align}\label{emeancell1}
\delta_0(p)=
-\gamma_2 (1-p)+(\gamma_0+\gamma_2)(1-p)^2-\gamma_0\sum^\infty_{m=3}p_{0,2}(m)(1-p)^m.
\end{align}
For edge percolation on $X$, 
\begin{align}\label{emeanedge}
\delta_0(p)=
\gamma_0 -(\gamma_0+\gamma_2)p-\gamma_0\sum^\infty_{m=3}p_{0,1}(m)(1-p)^m
+\gamma_2\sum^\infty_{m=3}p_{2,1}(m)p^m.
\end{align}
For vertex percolation on $X$, 
\begin{align}\label{emeanvertex}
\delta_0(p)=
\gamma_0 p-(\gamma_0+\gamma_2) p^2+\gamma_2\sum^\infty_{m=3}p_{2,0}(m)p^m.
\end{align}
\end{corollary}

\begin{example}\label{eArch1} \rm
An {\em Archimedean lattice} $A$
is a tessellation of the plane based on a finite number
of regular polygons such that all vertices are equivalent in
a graph-theoretical sense, see e.g.\ \cite[pp.154]{BR06b}.
It can be conveniently denoted by $(n_1,\ldots,n_z)$,
where $z$ is the degree of the vertices (called
{\em coordination number}) and $n_1,\ldots,n_z$ are the number
of edges of the polygons surrounding a vertex.
An Archimedean lattice can be made stationary
by putting $X:=A+\xi$ where $\xi$ is uniformly distributed
on a \emph{fundamental domain}, i.e.\ a connected set such that the tessellation 
can be generated by translations of the fundamental domain.
It is then easy to see that $p_{0,2}(m)=\I\{m=z\}$ and
$$
\gamma_2 p_{2,0}(m)=\gamma_0\sum^z_{k=1}\I\{n_k=m\}\frac{1}{n_k}.
$$
Because the Archimedean lattices are planar, we obtain $p_{0,1}(m)=p_{0,2}(m)$ and 
$p_{2,1}(m)=p_{2,0}(m)$. Further, we have $n_{0,1}=z$ and with \eqref{vertices} 
it follows $\gamma_2=\gamma_0\frac{z-2}{2}$ and 
$\gamma_1=\gamma_0+\gamma_2=\gamma_0 z/2$.
Now, we obtain for cell percolation
$$ 
\delta_0(p) = -\gamma_0 \frac{z}{2} p(1-p) + \gamma_0(1-p) - \gamma_0(1-p)^z, 
$$
for edge percolation
$$ 
\delta_0(p)=
\gamma_0 - \gamma_0 \frac{z}{2}p - \gamma_0 (1-p)^z 
+ \sum_{m=3}^\infty \gamma_0p^m \sum_{k=1}^z \I\{n_k=m\} \frac{1}{n_k}, $$
and for vertex percolation
$$ 
\delta_0(p) = \gamma_0p - \gamma_0 \frac{z}{2}p^2 
+ \sum_{m=3}^\infty \gamma_0p^m \sum_{k=1}^z \I\{n_k=m\} \frac{1}{n_k}. 
$$
\end{example}

\begin{example}\label{e3.4}\rm For cell percolation on a planar and normal
tessellation $X$, 
\begin{align}\label{meannormal2}
\delta_0(p)=\gamma_2 p(1-p)(1-2p),\quad p\in[0,1].
\end{align}
In particular, $\delta_0(1-p)=-\delta_0(p)$ and $\delta_0(1/2)=0$.
\end{example}

If $X$ is a {\em line tessellation} (see \cite{SW}) then $p_{0,2}(m)=0$ for $m\ne 4$.

\begin{corollary}\label{c2} For cell percolation on
a line tessellation,
\begin{align}\label{meanline}
\delta_0(p)=\gamma_2 p(1-p)(p^2-3p+1),
\quad p\in[0,1].
\end{align}
\end{corollary}

\section{Second order properties of face percolation}
\label{secondnormal}

In this section we consider $n$-percolation $\Psi$ on
a face-to-face tessellation $X$ for fixed $n \in \{0,\ldots,d\}$. We are interested in the
limits
\begin{align}\label{sigmaij}
\sigma_{i,j}(p):=
\lim_{t\to\infty} V_d(W_t)^{-1} \CV(V_i(Z\cap W_t),V_j(Z\cap W_t))
\end{align}
for $i,j\in\{0,\ldots,d\}$, where $W_t:=t^{1/d}W$ and $W\in\cP^d$ is a fixed 
polytope with volume one that contains the origin in its interior. 
Note that this definition depends on $W$.
Our aim is to establish a set of assumptions 
guaranteeing that these asymptotic covariances exist.
It is not hard to see that the result must be
a polynomial in the percolation parameter $p$.
The coefficients, however, are complicated, and are determined
by the global fluctuation of the intrinsic volumes within
the face processes $X_0,\dots,X_d$ as well as by the local
geometry of $X$, which is independent of $W$.

For $k\in \{0,\ldots,d\}$ and $r \in \N$ we define a polynomial
$f_n^k(r,\cdot)$ on $[0,1]$ by
\begin{align*}
f_n^{k}(r,p) :=
\I\{k<n\} (1-(1-p)^{r}) + \I\{k \ge n\} p^{r}.
\end{align*}
We need to assume the existence of the limits
\begin{align}\label{rhoij}
\rho^{k,l}_{i,j}(p):=
\lim_{t\to\infty} V_d(W_t)^{-1}
\CV\bigg(&\int V_i(F(x)\cap W_t) f_n^k(|\cS_n(x)|,p)\,\eta^{(k)}(\mathrm{d}x), \notag\\
&\int V_j(F(x)\cap W_t) f^l_n(|\cS_n(x)|,p) \,\eta^{(l)}(\mathrm{d}x)\bigg),
\end{align}
for all $i,j,k,l\in\{0,\dots,d\}$ and $p\in[0,1]$. Again these limits depend on $W$.


Further, we need to assume that
\begin{align}\label{Varnull}
  \lim_{t \rightarrow \infty} t^{-1/2} \sum_{r=1}^\infty 
\sqrt{\BV\bigg(\int (-1)^{\dim(F(x) \cap U_t)} V_i(F(x) \cap U_t) \I\{|\cS_n(x)|=r\} 
\, \eta^{(k)}(\mathrm{d}x)\bigg)}=0
\end{align}
for $i,k \in \{0,\ldots,d\}$ and $U \in \cF(W)$ with $\dim(U)<d$, 
where as before $U_t := t^{1/d}U$. We will use \eqref{Varnull} to control
the boundary term $V_i(Z\cap \partial W_t)$.
 
To describe the local neighbourhood of a point $x\in\R^d$ we take $l\in\{0,\dots,d\}$, 
$m\in\N$, and define $\cS^{m}_{l}(x)$
as the system of all $l$-dimensional faces sharing $m$ neighbouring $n$-faces with
the face $F(x)$, that is
\begin{align*}
\cS^{m}_{l}(x):=\{G\in X_l:|\cS_n(x) \cap \cS_n(s(G))|=m\}
\end{align*}
and
$$ \cS_l^{m,s}(x) := \{G \in \cS_l^m(x): |\cS_n(s(G))|=s\}, \quad s \in \N, $$
is the system of all $l$-faces in $\cS_l^m(x)$ that have $s$ neighbouring $n$-faces. 
We will assume that
\begin{align}\label{5.51}
\sum_{i,j,k=0}^d \gamma_n \BE_n^0[V_i(\cS_k(0))^2 V_j(\cS_d(0))] < \infty,
\end{align}
where, for any finite $\mathcal{S}\subset\cK^d$,
$$
V_i(\mathcal{S}):=\sum_{G\in \mathcal{S}}V_i(G)
$$
is the total $i$-th intrinsic volume of the members of $\cS$. For $j=0$, 
\eqref{5.51} implies \eqref{locfinite2} because \eqref{2.17} yields for 
$i,k \in \{0,\ldots,d\}$
\begin{align*}
\gamma_k \BE_k^0[V_i(F(0))^2] &\leq  \gamma_k \BE_k^0[V_i(F(0))^2 V_0(\cS_n(0))] 
\le \gamma_n \BE_n^0[V_i(\cS_k(0))^2 V_0(F(0))]\\
&\le \gamma_n \BE_n^0[V_i(\cS_k(0))^2 V_0(\cS_d(0))].
\end{align*}
In the following theorem we use the polynomial
\begin{align*}
g_{n,m}^{k,l,r,s}(p) := \, &\I\{k,l < n\} (1-p)^{r+s-m} (1-(1-p)^{m}) 
+ \I\{k \geq n, l < n\} p^{r} (1-p)^{s}\\
&+ \I\{k < n, l \geq n\} (1-p)^{r} p^{s} + \I\{l,k \geq n\} p^{r+s-m} (1-p^{m}),
\end{align*}
where $k,l \in \{0,\ldots,d\}$, $r,s \in \N$ and $m \in \{1,\ldots, \min(r,s)\}$.

\begin{theorem}\label{t5.1} Let $i,j\in\{0,\dots,d\}$. Assume that \eqref{Varnull} 
and \eqref{5.51} hold and the limits \eqref{rhoij} exist.
Then the limits \eqref{sigmaij} exist and are given by
\begin{align}\label{5.7}
\sigma_{i,j}(p) = \sum_{k=i}^d \sum_{l=j}^d &\,(-1)^{i+j+k+l}  
\bigg(\rho_{i,j}^{k,l}(p) \\ \notag
&+ \sum_{r,s=1}^\infty\sum_{m=1}^{\min(r,s)} g_{n,m}^{k,l,r,s}(p) \, \gamma_k \BE_k^0 \big[V_i(F(0)) 
V_j(\cS_l^{m,s}(0)) \I\{|\cS_n(0)|=r\} \big]\bigg).
\end{align}
\end{theorem}


{\em Proof.} First we show that the boundary term $V_i(Z \cap \partial W_t)$ 
is negligible, that is
\begin{align} \label{VorVarRand}
\lim_{t \to \infty} t^{-1} \BV(V_i(Z \cap \partial W_t))=0.
\end{align}
By \eqref{FormelRand} and 
the Cauchy-Schwarz inequality, it is enough to show that
\begin{align} \label{5.6a}
\lim_{t \to \infty} \frac{1}{t} \BV\bigg(\sum_{F \in X_k^1} (-1)^{i+\dim(F \cap U_t)} 
V_i(F \cap U_t)\bigg)=0
\end{align}
for $k \in \{0,\ldots,d\}$ and $U \in \cF(W)$ with $\dim(U)<d$. For 
the second moment of this functional we have
\begin{align}\label{e5.8}
&\BE\bigg[\bigg(\sum_{F \in X_k^1} (-1)^{i+\dim(F \cap U_t)} V_i(F \cap U_t)\bigg)^2\bigg]\\
&=
\begin{aligned}[t]
\sum_{r,s=1}^\infty \sum_{m=0}^{\min(r,s)} &\BE \iint (-1)^{\dim(F(x) \cap U_t) + \dim(F(y) \cap U_t)} 
V_i(F(x) \cap U_t) V_i(F(y) \cap U_t)\\ \notag
&\times \I\{|\cS_n(x)|=r, F(y) \in \cS_k^{m,s}(x)\} \, \I\{F(x),F(y) \in X_k^1\} \, 
\eta^{(k)}(\mathrm{d}x) \, \eta^{(k)}(\mathrm{d}y).
\end{aligned}
\end{align}
Take $x,y\in\eta^{(k)}$ and assume that $|\cS_n(x)|=r$ and 
$F(y) \in \cS_k^{m,s}(x)$ for $r,s\ge 1$ and $m\in\{0,\dots,\min(r,s)\}$.
The conditional probability of 
$\{F(x)\in X_k^1\}\cap \{F(y)\in X_k^1\}$ given $X$ is given by 
$(1-(1-p)^m+(1-p)^m (1-(1-p)^{r-m}) (1-(1-p)^{s-m}))$ in the case $k<n$ and
by $p^{r+s-m}$ in the case $k\ge n$. Therefore \eqref{e5.8} equals
\begin{align*}
&\begin{aligned}[t]
\sum_{r,s=1}^\infty \sum_{m=0}^{\min(r,s)} &\bigg( \I\{k<n\} (1-(1-p)^m 
+ (1-p)^m (1-(1-p)^{r-m}) (1-(1-p)^{s-m}))\\
&\quad + \I\{k \geq n\} p^{r+s-m} \bigg)\\
&\begin{aligned}[t]
\times \BE\bigg[ \iint & (-1)^{\dim(F(x) \cap U_t) + \dim(F(y) \cap U_t)} V_i(F(x) \cap U_t) 
V_i(F(y) \cap U_t)\\
&\times \I\{|\cS_n(x)|=r, F(y) \in \cS_k^{m,s}(x)\} \, \eta^{(k)}(\mathrm{d}x) 
\, \eta^{(k)}(\mathrm{d}y) \bigg].
\end{aligned}
\end{aligned}
\end{align*}
Since $\I\{F(y) \in \cS_k^{0,s}(x)\} = 1- \sum_{m=1}^{\min(r,s)} \I\{F(y) \in \cS_k^{m,s}(x)\}$, 
we get with \eqref{EWViRand}
\begin{align*}
&\BV\bigg(\sum_{F \in X_k^1} (-1)^{i+\dim(F \cap U_t)} V_i(F \cap U_t) \bigg)\\
&=
\begin{aligned}[t]
\sum_{r,s=1}^\infty
&\begin{aligned}[t]
f_n^k(r,p) f_n^k(s,p)
\CV\bigg( &\int (-1)^{\dim(F(x) \cap U_t)} V_i(F(x) \cap U_t) \I\{|\cS_n(x)|=r\} \, 
\eta^{(k)}(\mathrm{d}x),\\
& \int (-1)^{\dim(F(x) \cap U_t)} V_i(F(x) \cap U_t) \I\{|\cS_n(x)|=s\} \, 
\eta^{(k)}(\mathrm{d}x) \bigg)
\end{aligned}
\end{aligned}\\
&\quad+
\begin{aligned}[t]
\sum_{r,s=1}^\infty \sum_{m=1}^{\min(r,s)} g_{n,m}^{k,k,r,s}(p) \, 
\BE \iint &(-1)^{\dim(F(x) \cap U_t) + \dim(F(y) \cap U_t)} V_i(F(x) \cap U_t) V_i(F(y) \cap U_t)\\
&\times \I\{|\cS_n(x)|=r, F(y) \in \cS_k^{m,s}(x)\} \, \eta^{(k)}(\mathrm{d}x) \, 
\eta^{(k)}(\mathrm{d}y).
\end{aligned}
\end{align*}
Using the Cauchy-Schwarz inequality, $|f_n^{k}(r,p)|\leq 1$
and assumption \eqref{Varnull} we see that the first summand tends to 
zero after dividing by $t$. Now we consider the 
second summand. With $|g_{n,m}^{k,k,r,s}(p)| \leq 1$, the monotonicity of the 
intrinsic volumes, the monotone convergence theorem and the refined Campbell 
theorem \eqref{refCampbell}, we have
\begin{align*}
&\begin{aligned}[t]
\bigg|\sum_{r,s=1}^\infty \sum_{m=1}^{\min(r,s)} g_{n,m}^{k,k,r,s}(p) \frac{1}{t} 
\BE \iint &(-1)^{\dim(F(x) \cap U_t) + \dim(F(y) \cap U_t)} V_i(F(x) \cap U_t) V_i(F(y) \cap U_t)\\
&\times \I\{|\cS_n(x)|=r, F(y) \in \cS_k^{m,s}(x)\} \, \eta^{(k)}(\mathrm{d}x) \, 
\eta^{(k)}(\mathrm{d}y) \bigg|
\end{aligned}\\
&\leq
\begin{aligned}[t]
\sum_{r,s=1}^\infty \sum_{m=1}^{\min(r,s)} \frac{1}{t} \BE \iint &V_i(F(x)) \, V_i(F(y)) \, 
\mathbf{1}\{F(x) \cap U_t \neq \emptyset, F(y) \cap U_t \neq \emptyset\}\\
&\times \I\{|\cS_n(x)|=r, F(y) \in \cS_k^{m,s}(x)\} \, \eta^{(k)}(\mathrm{d}x) \, 
\eta^{(k)}(\mathrm{d}y)
\end{aligned}\\
&\leq
\begin{aligned}[t]
\frac{1}{t} \BE \iiint &V_i(F(x)) V_i(F(y)) \I\{F(x) \in \cS_k(z), F(y) \in \cS_k(z)\}\\
&\times \I \bigg\{\bigcup_{C \in \cS_d(z)}C \cap U_t \neq \emptyset \bigg\} \, 
\eta^{(k)}(\mathrm{d}x) \, \eta^{(k)}(\mathrm{d}y) \, \eta^{(n)}(\mathrm{d}z)
\end{aligned}\\
&\leq \frac{1}{t} \BE \int V_i(\cS_k(z))^2 
\I\bigg\{\bigcup_{C \in \cS_d(z)} C \cap \partial W_t \neq \emptyset\bigg\} \, 
\eta^{(n)}(\mathrm{d}z)\\
&= \gamma_n \BE_n^0 \bigg[V_i(\cS_k(0))^2 
\lambda_d\bigg(\partial W - t^{-1/d} \bigcup_{C \in \cS_d(0)} C\bigg)\bigg]\\
&\leq \gamma_n \BE_n^0 \bigg[V_i(\cS_k(0))^2 \sum_{C \in \cS_d(0)} 
\lambda_d (\partial W - t^{-1/d} C)\bigg].
\end{align*}
By the dominated 
convergence theorem this tends to zero as $t \rightarrow \infty$. Indeed,
for $t\ge 1$ 
$$
V_i(\cS_k(0))^2 \sum_{C \in \cS_d(0)} \lambda_d (W - C)
$$
is a dominating random variable.
The Steiner formula and  assumption \eqref{5.51} imply 
the required integrability of this random variable.
Hence \eqref{VorVarRand} follows.

Using \eqref{VorVarRand} and the representation \eqref{FormelVi}
together with the Cauchy-Schwarz inequality, we obtain 
\begin{align}\notag
\sigma_{i,j}(p) &= \lim_{t \to \infty} \frac{1}{t} 
\CV(V_i(Z \cap \mathrm{int}(W_t)), V_j(Z \cap \mathrm{int}(W_t)))\\
\label{e5.11}
&= \sum_{k,l=0}^d (-1)^{i+j+k+l} \lim_{t \rightarrow \infty} \frac{1}{t} 
\CV\bigg(\sum_{F \in X_k^1} V_i(F \cap W_t), \sum_{G \in X_l^1} V_j(G \cap W_t)\bigg).
\end{align}
The definition of $n$-percolation yields
\begin{align} \label{5.6}
\BE\bigg[\sum_{F \in X_k^1} V_i(F \cap W_t)\bigg] = \sum_{r=1}^\infty \, 
&( \I\{k<n\} (1-(1-p)^r) + \I\{k \geq n\} p^r) \notag\\
&\times \BE \int V_i(F(x) \cap W_t) \I\{|\cS_n(x)|=r\} \, \eta^{(k)}(\mathrm{d}x).
\end{align}
As \eqref{5.51} implies \eqref{locfinite2}, all occuring expectations are finite. 
To treat the mixed expectations \eqref{e5.11} we argue similarly
as in the first part of the proof. We have
\begin{align} \label{EW5.7}
&\BE\bigg[\sum_{F \in X_k^1} \sum_{G \in X_l^1} V_i(F \cap W_t) V_j(G \cap W_t)\bigg] \\
&=
\begin{aligned}[t]
\sum_{r,s=1}^\infty \sum_{m=0}^{\min(r,s)}
&\begin{aligned}[t]
\bigg( &\I\{k,l < n\} (1- (1-p)^{r} - (1-p)^{s} + (1-p)^{r+s-m})\\
&+ \I\{k \geq n,l < n\} (p^{r} - p^{r}(1-p)^{s} \I\{m=0\})\\
&+ \I\{k < n,l \geq n\} (p^{s} - p^{s}(1-p)^{r} \I\{m=0\})+ \I\{k,l \geq n\} p^{r+s-m}\bigg)
\end{aligned}\\
&\begin{aligned}[t]
\times \BE \iint &V_i(F(x) \cap W_t) V_j(F(y) \cap W_t) \I\{F(y) \in \cS_l^{m,s}(x)\}\\
&\times \I\{|\cS_n(x)|=r\} \, \eta^{(k)}(\mathrm{d}x) \, \eta^{(l)}(\mathrm{d}y).
\end{aligned}
\end{aligned} \notag
\end{align}
Since $\I\{F(y) \in \cS_l^{0,s}(x)\} = 1- \sum_{m=1}^{\min(r,s)} \I\{F(y) \in \cS_l^{m,s}(x)\}$, 
we get from \eqref{5.6} and \eqref{EW5.7} 
\begin{align}\label{5.12}
&\CV\bigg(\sum_{F \in X_k^1} V_i(F \cap W_t), \sum_{G \in X_l^1} V_j(G \cap W_t)\bigg) \notag\\
&=
\begin{aligned}[t]
\sum_{r,s=1}^\infty &f_n^k(r,p) f_n^l(s,p)\\
&\times \CV\bigg(\sum_{F \in X_k} V_i(F \cap W_t) 
\I\{|\cS_n(s(F))|=r\}, \sum_{G \in X_l} V_j(G \cap W_t) \I\{|\cS_n(s(G))|=s\}\bigg)
\end{aligned} \notag\\
&\quad+ \sum_{r,s=1}^\infty \sum_{m=1}^{\min(r,s)}
\begin{aligned}[t]
g_{n,m}^{k,l,r,s}(p) \, \BE \iint &V_i(F(x) \cap W_t) V_j(F(y) \cap W_t) 
\I\{F(y) \in \cS_l^{m,s}(x)\}\\
&\times \I\{|\cS_n(x)|=r\} \, \eta^{(k)}(\mathrm{d}x) \, \eta^{(l)}(\mathrm{d}y).
\end{aligned}
\end{align}
By assumption \eqref{rhoij} the first summand tends to $\rho_{i,j}^{k,l}(p)$ 
after dividing by $t$.  
Using $|g_{n,m}^{k,l,r,s}(p)| \leq 1$, the monotone convergence theorem and the 
refined Campbell theorem we get
\begin{align*}
&\begin{aligned}[t]
\bigg| \sum_{r,s=1}^\infty \sum_{m=1}^{\min(r,s)} g_{n,m}^{k,l,r,s}(p) \frac{1}{t} \bigg(&
\begin{aligned}[t]
\BE \iint &V_i(F(x) \cap W_t) V_j(F(y) \cap W_t) \I\{F(y) \in \cS_l^{m,s}(x)\}\\
&\times \I\{|\cS_n(x)|=r\} \, \eta^{(k)}(\mathrm{d}x) \, \eta^{(l)}(\mathrm{d}y)
\end{aligned}\\
&-
\begin{aligned}[t]
\BE \iint &\I\{x \in W_t\} V_i(F(x)) V_j(F(y)) \I\{F(y) \in \cS_l^{m,s}(x)\}\\
&\times \I\{|\cS_n(x)|=r\} \, \eta^{(k)}(\mathrm{d}x) \, \eta^{(l)}(\mathrm{d}y) \bigg) \bigg|
\end{aligned}
\end{aligned}\\
&\leq
\begin{aligned}[t]
\frac{1}{t} \BE \iiint &|V_i(F(x) \cap W_t) V_j(F(y) \cap W_t) - \I\{x \in W_t\} 
V_i(F(x)) V_j(F(y))|\\
&\times \I\{F(x) \in \cS_k(z), F(y) \in \cS_l(z)\} \, \eta^{(k)}(\mathrm{d}x) \, 
\eta^{(l)}(\mathrm{d}y) \, \eta^{(n)}(\mathrm{d}z).
\end{aligned}
\end{align*}
This can be further bounded by
$$
\gamma_n \BE_n^0\bigg[
\lambda_d\bigg(\partial W - t^{-1/d} \bigcup_{C \in \cS_d(0)} C \bigg) V_i(\cS_k(0)) 
V_j(\cS_l(0)) \bigg].
$$
Because of the monotonicity of the Lebesgue measure, the Steiner formula, 
the Cauchy-Schwarz inequality and assumption \eqref{5.51} this converges to 0 as 
$t \rightarrow \infty$ by the dominated convergence theorem. Hence, 
by applying the refined Campbell theorem the second summand of \eqref{5.12} 
converges after dividing by $t$ to the second summand of the assertion.\qed

\bigskip

For cell percolation on a normal tessellation $X$ we can find similarly to 
Proposition \ref{BezEW} a relation between $\sigma_{i,j}(p)$ and $\sigma_{i,j}(1-p)$.

\begin{proposition}\label{BezCov}
Consider cell percolation on a normal tessellation $X$. Assume that \eqref{Varnull} 
and \eqref{5.51} hold and the limits \eqref{rhoij} exist. Then we have for $p \in [0,1]$
\begin{align*}
\sigma_{i,j}(p) = (-1)^{i+j} \, \sigma_{i,j}(1-p), \quad i,j \in \{0,\ldots,d\}.
\end{align*}
\end{proposition}
From \eqref{Z(1-p)} and \eqref{Z(1-p)2} we obtain that
\begin{align*}
(-1)^{i+j} \, &\CV(V_i(Z_{1-p} \cap W_t), V_j(Z_{1-p} \cap W_t))\\
&= \CV(V_i(Z_p \cap \mathrm{int}(W_t)), V_j(Z_p \cap \mathrm{int}(W_t))).
\end{align*}
Now the assertion follows from \eqref{VorVarRand} and the Cauchy-Schwarz inequality.
\qed

\bigskip

Formula \eqref{5.7} simplifies for cell percolation on a normal tessellation. 
We define for $i,j,k,l \in \{0,\ldots,d\}$
\begin{align}\label{rhoohners}
\tau_{i,j}^{k,l} 
:=\lim_{t \to \infty} V_d(W_t)^{-1} \CV\bigg( \int V_i(F(x) \cap W_t) \, 
\eta^{(k)}(\mathrm{d}x), \int V_j(F(x) \cap W_t) \, \eta^{(l)}(\mathrm{d}x) \bigg).
\end{align}

\begin{proposition}\label{cellnormalvar}
Consider cell percolation on a normal tessellation $X$ and let $i,j \in \{0,\ldots,d\}$. 
Assume that \eqref{Varnull} and \eqref{5.51} hold and the limits \eqref{rhoohners} exist. 
Then we have for $p \in [0,1]$
\begin{align}\label{cellnormalkov}
\sigma_{i,j}(p) = \sum_{k=i}^d \sum_{l=j}^d (-1)^{k+l}\bigg(p^{2d-k-l+2} \, \tau_{i,j}^{k,l}
+ \sum_{m=1}^{d-\max(k,l)+1} &p^{2d-k-l-m+2} (1-p^m)\\
&\times \gamma_k \, \BE_k^0[V_i(F(0)) \, V_j(\cS_l^m(0))] \bigg).
\notag
\end{align}
\end{proposition}
\emph{Proof.} We use Theorem \ref{t5.1} and first note that the limits 
\eqref{rhoij} exist and are given by
$$ 
\rho_{i,j}^{k,l}(p) = (1-(1-p)^{d-k+1})(1-(1-p)^{d-l+1})\tau_{i,j}^{k,l}. 
$$
Since $n=d$, the normality implies almost surely that $|\cS_n(x)|=d-k+1$ 
for $x \in \eta^{(k)}$ and $k \in \{0,\ldots,d\}$. Moreover, it is
easy to check that
\begin{align*}
g^{k,l,r,s}_{n,m}(p)&=(1-p)^{2d-k-l-m+2}(1-(1-p)^m).
\end{align*}
By the inclusion-exclusion principle,
\begin{align*}
\sum_{k=i}^d (-1)^{d-k} \int V_i(F(x) \cap W_t) \, \eta^{(k)}(\mathrm{d}x) = V_i(W_t),
\end{align*}
so that $\sum^d_{k=i}(-1)^k \tau^{k,l}_{i,j}=\sum^d_{l=j}(-1)^k \tau^{k,l}_{i,j}=0$.
Therefore formula \eqref{5.7} yields
\begin{align*}
\sigma_{i,j}(p) = &\sum_{k=i}^d \sum_{l=j}^d (-1)^{i+j+k+l} \bigg( (1-p)^{2d-k-l+2} \, 
\tau_{i,j}^{k,l}\\
&+ \sum_{m=1}^{d-\max(k,l)+1} (1-p)^{2d-k-l-m+2} (1-(1-p)^m) 
\, \gamma_k \, \BE_k^0[V_i(F(0)) \, V_j(\cS_l^m(0))] \bigg).
\end{align*}
Combining this with Proposition \ref{BezCov} yields the assertion.
\qed

\section{On the covariance structure in the plane}\label{secvEuler}

In this section we consider cell percolation on a planar and normal tessellation
$X$ satisfying  \eqref{Varnull} and \eqref{5.51} and assume that the 
limits \eqref{rhoohners} exist.
We define 
the expected square of the number of vertices of a typical cell by
\begin{align}\label{mu2}
\mu_2:=\BE^0_2 f_0(F(0))^2.
\end{align}
Here $f_0(P)$ denotes the number of vertices of a polygon $P\subset\R^2$. 
Since $\BE_2^0 f_0(F(0))=6$ (see Remark \ref{r1}),
Jensen's inequality gives
\begin{align}\label{36}
\mu_2 \ge 36.
\end{align}


\begin{example}\rm \label{PVoronoi} Assume that 
$X$ is the Voronoi tessellation generated by a stationary
Poisson process. In this case integral expressions
for $\mu_2$ are available. Numerical integration
gives $\mu_2\approx 37.78$, see \cite{HeiMu08}.
Therefore $\alpha\approx 13.89$.
\end{example}

The following main result of this section expresses the asymptotic convariance
structure in terms of second order properties of the
typical cell and the typical edge.
Recall the definition \eqref{rhoohners} of $\tau_{1,1}^{2,2}$,
$\tau_{1,0}^{2,2}$, and $\tau_{0,0}^{2,2}$.

\begin{theorem}\label{Kovarianzstruktur}
Assume that the limits \eqref{rhoohners} exist and that \eqref{Varnull} and \eqref{5.51} 
are satisfied. Then the asymptotic covariance structure is given by
\begin{align*}
\sigma_{2,2}(p) &= p(1-p) \gamma_2 \BE_2^0 V_2(F(0))^2,\\
\sigma_{1,2}(p) &= p(1-p)(1-2p) \gamma_2 \BE_2^0 [V_2(F(0))V_1(F(0))],\\
\sigma_{0,2}(p) &= p(1-p) - p^2(1-p)^2 \gamma_2 \BE_2^0[V_2(F(0)) f_0(F(0))],\\
\sigma_{1,1}(p) &= p^2(1-p)^2 (\tau_{1,1}^{2,2} + \gamma_1 \BE_1^0V_1(F(0))^2) 
+ p(1-p)(1-2p)^2 \gamma_2 \BE_2^0V_1(F(0))^2,\\
\sigma_{0,1}(p) &= p^2(1-p)^2(1-2p) (\tau_{1,0}^{2,2} 
-\gamma_2 \BE_2^0[V_1(F(0)) f_0(F(0))])\\
&\quad+ p(1-p)(1-p-3p^2+2p^3) \gamma_2 \BE_2^0[V_1(F(0))],\\
\sigma_{0,0}(p) &= \gamma_2 \mu_2 p^3(1-p)^3 + \gamma_2 p(1-p)(1-9p-p^2+20p^3-10p^4)\\
&\quad+ \tau_{0,0}^{2,2} p^2(1-p)^2(1-2p)^2.
\end{align*}
\end{theorem}
{\em Proof.} The formulae for $\sigma_{2,2}$, $\sigma_{1,2}$ and $\sigma_{0,2}$ 
follow directly from Proposition \ref{cellnormalvar} by using 
$\gamma_2 \, \BE_2^0 V_2(F(0))=1$. 

To treat $\sigma_{1,1}$ we first recall that, for any convex body
$K\subset\R^2$,  $V_1(K)=\frac{1}{2}\mathcal{H}^1(\partial K)$
if $K$ has non-empty interior. Otherwise
$V_1(K)=\mathcal{H}^1(K)$. Here $\mathcal{H}^1$ denotes 
the one-dimensional Hausdorff measure on $\R^2$. It follows that
\begin{align}\label{Gleich}
\int V_1(F(x) \cap W_t) \, \eta^{(1)}(\mathrm{d}x) = \int V_1(F(x) \cap W_t) \, 
\eta^{(2)}(\mathrm{d}x) - \frac{1}{2}\mathcal{H}^1(\partial W_t),
\end{align}
and therefore
$$ 
\tau_{1,1}^{1,1} = \tau_{1,1}^{1,2} = \tau_{1,1}^{2,2}. 
$$
Proposition \ref{t2.1} yields after a straightforward calculation that
\begin{align*}
\gamma_1 \BE_1^0[V_1(F(0)) \, V_1(\cS_1^1(0))] &= 4 \gamma_2 \BE_2^0[V_1(F(0))^2] 
- 2\gamma_1 \BE_1^0[V_1(F(0))^2],\\
\gamma_1 \BE_1^0[V_1(F(0)) \, V_1(\cS_2^1(0))] &= 2 \gamma_2 \BE_2^0[V_1(F(0))^2].
\end{align*}
Inserting the above formula into Proposition \ref{cellnormalvar} 
yields the asserted formula for $\sigma_{1,1}$ after a simple calculation.

To deal with the remaining covariances $\sigma_{0,1}$ and 
$\sigma_{0,0}$, we define
\begin{align}\label{Defeps}
\varepsilon_t := \sum_{e \in X_1} V_0(e \cap \partial W_t), \quad t>0.
\end{align}
Taking $i=0$, $k=1$ and $r=2$ in \eqref{Varnull} we obtain
\begin{align}\label{Voreps}
\lim_{t \rightarrow \infty} t^{-1} \mathrm{Var}(\varepsilon_t)=0.
\end{align}
Euler's formula yields
\begin{align*}
(|X_{0,t}|+\varepsilon_t) + (|X_{2,t}|+1) = (|X_{1,t}|+\varepsilon_t)+2,
\end{align*}
where $X_{k,t}$ denotes the set of all $k$-faces that have non-empty intersection 
with $W_t$. Further, by normality we have
\begin{align*}
2(|X_{1,t}|+\varepsilon_t) = 3(|X_{0,t}|+\varepsilon_t).
\end{align*}
Combinig these two equations yields
\begin{align}\label{EulerFormel}
|X_{0,t}|=2|X_{2,t}|-\varepsilon_t-2,\quad |X_{1,t}|=3|X_{2,t}|-\varepsilon_t-3.
\end{align}
With these observations we can determine $\sigma_{0,1}$ and $\sigma_{0,0}$.

Using \eqref{EulerFormel} and assumption \eqref{Voreps} we obtain
\begin{align*}
\tau_{1,0}^{1,0}=\tau_{1,0}^{2,0}=2 \tau_{1,0}^{2,2},\quad \tau_{1,0}^{1,1}
=\tau_{1,0}^{2,1}=3 \tau_{1,0}^{2,2}, \quad \tau_{1,0}^{1,2}=\tau_{1,0}^{2,2}
\end{align*}
and with Proposition \ref{t2.1} we get
\begin{align*}
\gamma_1 \BE_1^0[V_1(F(0)) \, V_0(\cS_0^1(0))] 
&= 2 \gamma_2 \BE_2^0[V_1(F(0)) \, f_0(F(0))] - 4 \gamma_2 \BE_2^0[V_1(F(0))],\\
\gamma_1 \BE_1^0[V_1(F(0)) \, V_0(\cS_1^1(0))] 
&=2\gamma_2 \BE_2^0[V_1(F(0)) \, f_0(F(0))] - 2 \gamma_2 \BE_2^0[V_1(F(0))],\\
\gamma_1 \BE_1^0[V_1(F(0))] &= \gamma_2 \BE_2^0[V_1(F(0))].
\end{align*}
Together with Proposition \ref{cellnormalvar}, these observations yield 
the asserted formula for $\sigma_{0,1}$.

Next, we determine $\tau_{0,0}^{k,l}$. Again with \eqref{EulerFormel} 
and assumption \eqref{Voreps} we obtain
\begin{align*}
\tau_{0,0}^{0,0}=4 \tau_{0,0}^{2,2}, \quad \tau_{0,0}^{0,1}=6 \tau_{0,0}^{2,2}, 
\quad \tau_{0,0}^{0,2}=2 \tau_{0,0}^{2,2}, \quad \tau_{0,0}^{1,1}=9 \tau_{0,0}^{2,2}, 
\quad \tau_{0,0}^{1,2}=3 \tau_{0,0}^{2,2}.
\end{align*}

To determine the second summand of $\sigma_{0,0}$ define
\begin{align*}
f(k,l,m):= \gamma_k \BE_k^0|\cS_{l}^m(0)|
\end{align*}
for $k,l\in\{0,1,2\},m\in\{1,\ldots,3-\max(k,l)\}$. 
Using Proposition \ref{t2.1} together with normality yields
\begin{align*}
f(0,0,1)&=\gamma_0 \BE_0^0\bigg[\sum_{G \in \cS_2(0)} f_0(G)-9\bigg] 
= \gamma_2 \mu_2 -9 \gamma_0,\\
f(1,0,1)&= \gamma_1 \BE_1^0\bigg[\sum_{G \in \cS_2(0)} f_0(G)-4\bigg] 
= \gamma_2 \mu_2 -4 \gamma_1,\\
f(0,1,1)&= \gamma_0 \BE_0^0\bigg[\sum_{G \in \cS_2(0)} f_0(G) -6\bigg] 
= \gamma_2 \mu_2-6\gamma_0,\\
f(1,1,1)&= \gamma_1 \BE_1^0\bigg[\sum_{G \in \cS_2(0)} f_0(G) -2\bigg] 
= \gamma_2 \mu_2 -2 \gamma_1,\\
f(0,0,2)&=f(0,1,2)=f(2,0,1)=f(0,2,1)=f(2,1,1)=3\gamma_0,\\
f(1,2,1)&=f(1,0,2)=2\gamma_1,\\
f(0,0,3)&=\gamma_0, \quad f(1,1,2)=\gamma_1, \quad f(2,2,1)=\gamma_2.
\end{align*}
Using Proposition \ref{cellnormalvar} and the relations $\gamma_0=2\gamma_2$ 
and $\gamma_1=3\gamma_2$ yields the assertion.
\qed \bigskip

Obviously, 0 and 1 are zeroes of all covariances considered in Theorem 
\ref{Kovarianzstruktur}. We continue with a brief discussion of the maxima and 
minima.

\begin{corollary}
Let the assumptions of Theorem \ref{Kovarianzstruktur} be satisfied. 
Then the variance $\sigma_{2,2}$ has a global maximum at 1/2 and $\sigma_{1,2}$ 
has a global maximum at $\frac{1}{2}-\frac{1}{2\sqrt{3}}$ and a global minimum at 
$\frac{1}{2}+\frac{1}{2\sqrt{3}}$. The covariance $\sigma_{0,2}$ has a global minimum at 
1/2 and the variance $\sigma_{1,1}$ has a global maximum (minimum) at 1/2 if
$$ 
2 \gamma_2 \BE_2^0[V_1(F(0))^2] < (>) \tau_{1,1}^{2,2} + \gamma_1 \BE_1^0[V_1(F(0))^2]. 
$$
The variance $\sigma_{0,0}$ has a strict global maximum (minimum) at $1/2$ if
$$ 
\mu_2 >(<) \frac{86}{3}+\frac{4\tau_{0,0}^{2,2}}{3\gamma_2}. 
$$
\end{corollary}

\section{Poisson Voronoi percolation}\label{poissonvoronoi}

In this section we consider the Voronoi tessellation $X$ generated by
a stationary Poisson process $\eta$ in $\R^d$ with intensity $\gamma>0$.
For a formal definition we introduce the space $\bN$ of all locally
finite subsets $\mu$ of $\R^d$ whose convex hull coincides with $\R^d$
and whose points are in {\em general quadratic position}. The latter means that 
no $d+2$ points of $A$ lie on the boundary of some
ball and any $k\in\{2,\ldots,d+1\}$  points in $\mu$
do not lie in a $(k-2)$-dimensional affine subspace of $\R^d$. 
The {\em Voronoi cell} $C(\mu,x)$ of $x\in \mu\in\bN$ is the set of all
$y\in\R^d$ satisfying $|y-x|\le \min\{|y-z|:z\in \mu\}$. 
The system $\{C(\mu,x):x \in \mu\}$ of all Voronoi cells with respect to $\mu$ 
is called the \emph{Voronoi tessellation}.

We can assume without restriction of generality that $\eta(\omega)\in\bN$
for all $\omega\in\Omega$.
The {\em Poisson Voronoi tessellation} $X:=\{C(\eta,x):x\in\eta\}$ is then 
stationary, face-to-face, and
normal, see Theorems 10.2.2 and 10.2.3 in \cite{SW}.

For $x,y \in \R^d$, let us define $\eta^x:=\eta\cup\{x\}$ and $\eta^{x,y}:=\eta\cup\{x,y\}$.
To abbreviate our notation we define the random variables
\begin{align*}
V^{(k)}_i(x,p)&:= \sum_{F \in \cF_k(C(\eta^{x},x))} V_i(F) \, f_n^k(|\cS_n(s(F))|,p),\\
V_i^{(k)}(x,y,p) &:= \sum_{F \in \cF_k(C(\eta^{x,y},x))} V_i(F) \, f_n^k(|\cS_n(s(F))|,p)
\end{align*}
for $i,k\in\{0,\ldots,d\}$, $x,y \in \R^d$ and $p \in [0,1]$.

We now show, for fixed $n \in \{0,\ldots,d\}$, that the assumptions
of Theorem \ref{t5.1} are satisfied. Moreover, we obtain
a more explicit representation of the limits \eqref{rhoij}. 
In particular, these limits are independent of the observation window $W$.

\begin{theorem} \label{thm7.1}
The Poisson Voronoi tessellation satisfies \eqref{Varnull} and \eqref{5.51}.
Moreover, the limits \eqref{rhoij} exist and are given by
\begin{align}
(d-k+1)(d-l+1)&\rho^{k,l}_{i,j}(p) 
=\gamma \, \BE[V_i^{(k)}(0,p) V_j^{(l)}(0,p)] \label{rhoijkl}\\
&+ \gamma^2 \int \bigg( \BE[V_i^{(k)}(x,0,p) \, V_j^{(l)}(0,x,p)]- \BE V_i^{(k)}(0,p) \, 
\BE V_j^{(l)}(0,p) \bigg) \mathrm{d}x. \notag
\end{align}
\end{theorem}
{\em Proof.}
Assumption \eqref{5.51} is a consequence of the Cauchy-Schwarz inequality, 
Lemma \ref{Appendix2} and Lemma \ref{Umkugelradius2}.

Next, we will show that assumption \eqref{Varnull} holds. Because the Poisson 
Voronoi tessellation is normal, we have for $i,k \in \{0,\ldots,d\}$, $r \in \N$ 
and $U \in \cF(W)$ with $\dim(U)<d$
\begin{align}\label{normalumschr}
&\int (-1)^{i+\dim(F(x) \cap U_t)} V_i(F(x) \cap U_t) \I\{|\cS_n(x)|=r\} 
\, \eta^{(k)}(\mathrm{d}x) \notag\\
&= \frac{1}{d-k+1} \int \sum_{F \in \cF_k(C(\eta,x))} (-1)^{i+\dim(F \cap U_t)} 
V_i(F \cap U_t) \I\{|\cS_n(s(F))|=r\} \, \eta(\mathrm{d}x).
\end{align}
Thus, it will be enough to show that
\begin{align}\label{7.3a}
\lim_{t \rightarrow \infty} \sum_{r=1}^\infty \sqrt{\frac{1}{t} 
\BV \bigg( \int \sum_{F \in \cF_k(C(\eta,x))} (-1)^{i+\dim(F \cap U_t)} V_i(F \cap U_t) 
\I\{|\cS_n(s(F))|=r\} \, \eta(\mathrm{d}x) \bigg)} =0.
\end{align}
To abbreviate the notation we define
\begin{align*}
h_r(\mu,B) := \int \sum_{F \in \cF_k(C(\mu,x))} (-1)^{i+\dim(F \cap B)} V_i(F \cap B) 
\I\{|\cS_n(s(F))|=r\} \, \mu(\mathrm{d}x)
\end{align*}
for $r \in \N$, $\mu \in \mathbf{N}$ and Borel sets $B \subset \R^d$. 
We use the Poincar\'{e} inequality, see \cite{Wu00}, to get
\begin{align}\label{poincare}
\frac{1}{t} \, \BV( h_r(\eta,U_t)) &\leq \frac{\gamma}{t} \, 
\BE \int (h_r(\eta^x,U_t) - h_r(\eta,U_t))^2 \, \mathrm{d}x.
\end{align}
Now we will determine an upper bound for $|h_r(\eta^x,U_t)-h_r(\eta,U_t)|$. 
Therefore, define the \emph{neighbourhood} of a point $x \in \mu$ with 
respect to $\mu\in \mathbf{N}$ by
\begin{align}\label{Nachbarschaft}
N(\mu,x):=\{y \in \mu \setminus \{x\}:C(\mu,x)\cap C(\mu,y)\neq \emptyset\}
\end{align}
and the \emph{neighbourhood of second order} of a point $x \in \mu$ with 
respect to $\mu$ by
\begin{align} \label{Nachbarschaft2}
N_2(\mu,x):=\{y \in \mu: \exists z \in N(\mu,x) \text{ with } y \in N(\mu,z)\}.
\end{align}
For $k \leq n$ a $k$-face is almost surely contained in $\binom{d-k+1}{n-k}$ 
$n$-faces since $X$ is normal and so, we have almost surely
$$ 
h_r(\eta,U_t) = h_{r}(\eta,U_t) \, \mathbf{1}\bigg\{r=\binom{d-k+1}{n-k}\bigg\}. 
$$
For $k>n$, the normality implies that each $n$-face of $F \in \cF_k(C(\eta,x))$ 
is contained in $C(\eta,x)$ and $d-n$ neighbouring cells of $x$. 
So, if $|\cS_n(s(F))|=r$, there must be at least $r$ ways to choose $d-n$ 
cells from the neighbouring cells of $x$ and this implies $|N(\eta,x)| \geq r^{1/(d-n)}$. 
Because the addition of a point $x \in \R^d$ to $\eta$ just changes the 
cells of the points $y \in \eta$ with $y \in N(\eta^x,x)$ we get with the 
previous observations
\begin{align*}
&|h_r(\eta^x,U_t) - h_r(\eta,U_t)|\\
&=
\begin{aligned}[t]
\bigg| &
\begin{aligned}[t]
\sum_{F \in \cF_k(C(\eta^x,x))} &(-1)^{i+\dim(F \cap U_t)}V_i(F\cap U_t)\I\{|\cS_n(\eta^x,s(F))|=r\}\\
& \times \bigg(\I\bigg\{k\le n,r=\binom{d-k+1}{n-k}\bigg\}+\I\{k>n,|N(\eta^x,x)| 
\ge r^{1/(d-n)}\} \bigg)
\end{aligned}\\
&+
\begin{aligned}[t]
\sum_{y \in N(\eta^x,x)} &\sum_{F \in \cF_k(C(\eta^x,y))} (-1)^{i+\dim(F \cap U_t)} 
V_i(F \cap U_t) \I\{|\cS_n(\eta^x,s(F))|=r\}\\
& \times \bigg(\I\bigg\{k \leq n, r = \binom{d-k+1}{n-k}\bigg\} 
+ \I\{k>n,|N(\eta^x,y)| \geq r^{1/(d-n)}\} \bigg)
\end{aligned}\\
&-
\begin{aligned}[t]
\sum_{y \in N(\eta^x,x)} &\sum_{F \in \cF_k(C(\eta,y))} (-1)^{i+\dim(F \cap U_t)} 
V_i(F \cap U_t) \I\{|\cS_n(s(F))|=r\}\\
& \times \bigg(\I\bigg\{k \leq n, r = \binom{d-k+1}{n-k}\bigg\} 
+ \I\{k>n,|N(\eta,y)| \ge r^{1/(d-n)}\} \bigg) \bigg|.
\end{aligned}
\end{aligned}
\end{align*}
In view of $\max(|N(\eta^x,x)|, |N(\eta^x,y)|, |N(\eta,y)|) \leq |N_2(\eta^x,x)|$ 
for $y \in N(\eta^x,x)$ and the monotonicity of the intrinsic volumes we have
\begin{align*}
&|h_r(\eta^x,U_t) - h_r(\eta,U_t)|\\
&\leq \sum_{F \in \cF_k(C(\eta^x,x))} V_i(F) \I\{F \cap U_t \neq \emptyset\}\\
&\quad+ \sum_{y \in N(\eta^x,x)} \bigg( \sum_{F \in \cF_k(C(\eta^x,y))} V_i(F) 
\I\{F \cap U_t \neq \emptyset\} + \sum_{F \in \cF_k(C(\eta,y))} V_i(F) 
\I\{F \cap U_t \neq \emptyset\} \bigg)\\
&\quad \times \bigg(\I\bigg\{k \leq n, r = \binom{d-k+1}{n-k}\bigg\} 
+ \I\{k>n,|N_2(\eta^x,x)| \geq r^{1/(d-n)}\}\bigg)\\
&\le \bigg( V_i(\cF_k(C(\eta^x,x))) + \sum_{y \in N(\eta^x,x)} (V_i(\cF_k(C(\eta^x,y))) 
+ V_i(\cF_k(C(\eta,y)))) \bigg)\\
&\quad \times \I \bigg\{\bigcup_{y \in N(\eta^x,x) \cup \{x\}} C(\eta^x,y) \cap U_t 
\ne \emptyset \bigg\}\\
&\quad \times \bigg(\I\bigg\{k \leq n, r = \binom{d-k+1}{n-k}\bigg\} 
+\I\{k>n,|N_2(\eta^x,x)| \geq r^{1/(d-n)}\}\bigg).
\end{align*}
To abbreviate the notation we define for $x \in \mu \in \mathbf{N}$ and fixed 
$k \in \{0,\ldots,d\}$
$$ 
f(\mu,x) := \bigg( V_i(\cF_k(C(\mu,x))) + \sum_{y \in N(\mu,x)} (V_i(\cF_k(C(\mu,y))) + 
V_i(\cF_k(C(\mu-\delta_x,y)))) \bigg). 
$$
Note that all moments of $f(\eta^0,0)$ exist by Lemma \ref{Nachbarzahl} 
and Lemma \ref{Umkugelradius2}.
Using \eqref{poincare}, the translation covariance of $C(\cdot,\cdot)$ and $\cF_k$, 
the stationarity of $\eta$ and the translation invariance of the intrinsic 
volumes and the number of neighbours yields
\begin{align*}
&\frac{1}{t} \BV(h_r(\eta,U_t))\\
&\leq
\begin{aligned}[t]
\frac{\gamma}{t} \, \BE \int &\I\bigg\{\bigcup_{y \in N(\eta^x,x)\cup\{x\}}
C(\eta^x,y) \cap U_t \neq \emptyset\bigg\} f(\eta^x,x)^2\\
&\times \bigg(\I\bigg\{k \leq n,r= \binom{d-k+1}{n-k}\bigg\} 
+ \I\{k>n,|N_2(\eta^x,x)| \geq r^{1/(d-n)}\}\bigg) \, \mathrm{d}x
\end{aligned}\\
&= \frac{\gamma}{t} \BE \int \I\bigg\{\bigg(\bigcup_{y \in N(\eta^0,0) \cup \{0\}}C(\eta^0,y)+x\bigg) 
\cap U_t \neq \emptyset\bigg\} \, f(\eta^0,0)^2\\
&\quad\times \bigg(\I\bigg\{k \leq n,r= \binom{d-k+1}{n-k}\bigg\} 
+ \I\{k>n,|N_2(\eta^0,0)| \geq r^{1/(d-n)}\}\bigg) \, \mathrm{d}x\\
&= \I\bigg\{k \leq n, r= \binom{d-k+1}{n-k}\bigg\} 
\, \gamma \BE \bigg[ V_d\bigg(U- t^{-1/d} \bigcup_{y \in N(\eta^0,0) \cup \{0\}}C(\eta^0,y)\bigg) 
\, f(\eta^0,0)^2 \bigg]\\
&\quad+ \I\{k>n\} \, \gamma \BE \bigg[ V_d\bigg(U- t^{-1/d} \hspace*{-1mm} \bigcup_{y \in N(\eta^0,0) \cup \{0\}} \hspace*{-1mm} C(\eta^0,y)\bigg) \, f(\eta^0,0)^2 \, \I\{|N_2(\eta^0,0)| \geq r^{1/(d-n)}\} \bigg].
\end{align*}
Applying the Cauchy-Schwarz inequality to the second summand, we get
\begin{align*}
&\frac{1}{t} \BV(h_r(\eta,U_t))\\
&\leq \I\bigg\{k \leq n, r= \binom{d-k+1}{n-k}\bigg\} \, 
\gamma \BE \bigg[ V_d\bigg(U- t^{-1/d} \bigcup_{y \in N(\eta^0,0) \cup \{0\}}C(\eta^0,y)\bigg) 
\, f(\eta^0,0)^2 \bigg]\\
&\quad+ \I\{k>n\} \, \gamma \BE \bigg[ V_d\bigg(U- t^{-1/d} 
\bigcup_{y \in N(\eta^0,0) \cup \{0\}}C(\eta^0,y)\bigg)^2 \, f(\eta^0,0)^4\bigg]^{1/2}\\
&\qquad \times \BP(|N_2(\eta^0,0)| \geq r^{1/(d-n)})^{1/2}.
\end{align*}
The expectations in both summands converge to 0 as $t \rightarrow \infty$ 
by the dominated convergence theorem and because of Lemma \ref{Nachbarzahl}, 
Lemma \ref{Umkugelradius2} and the Steiner formula a dominating function is given by
$$ 
V_d\bigg(U- \bigcup_{y \in N(\eta^0,0) \cup \{0\}}C(\eta^0,y)\bigg) \, f(\eta^0,0)^2 
$$
resp.\ its square. Hence, we get for $k \leq n$
\begin{align*}
\lim_{t \to \infty} \sum_{r=1}^\infty \sqrt{t^{-1} \BV(h_r(\eta,U_t))} 
&= \lim_{t \to \infty} \sqrt{t^{-1} \BV(h_{\binom{d-k+1}{n-k}}(\eta,U_t))} =0
\end{align*}
and for $k>n$
\begin{align*}
\lim_{t \to \infty} \sum_{r=1}^\infty \sqrt{t^{-1} \BV(h_r(\eta,U_t))} \leq \lim_{t \to \infty} \, 
&\gamma \BE \bigg[ V_d\bigg(V- t^{-1/d} \bigcup_{y \in N(\eta^0,0) \cup \{0\}}C(\eta^0,y)\bigg)^2 \, 
f(\eta^0,0)^4\bigg]^{1/4}\\
&\times \sum_{r=1}^\infty \BP(|N_2(\eta^0,0)| \geq r^{1/(d-n)})^{1/4}.
\end{align*}
Since the neighbourhood of second order of the typical point has an 
exponentially decreasing tail, see \eqref{Nach2tail}, this is zero, too. 
This proves \eqref{Varnull}.

We will prove \eqref{rhoijkl} in two steps. In the first, we consider 
an asymptotic covariance that is similar to $\rho_{i,j}^{k,l}(p)$ but easier to determine, 
i.e.\ we will show that
\begin{align}\label{firststep}
&\begin{aligned}[t]
\lim_{t \to \infty} \frac{1}{t} \CV\bigg(&\int V_i^{(k)}(x,p) \, \mathbf{1}\{x \in W_t\} \, 
\eta(\mathrm{d}x),\int V_j^{(l)}(x,p) \, \mathbf{1}\{x \in W_t\} \, 
\eta(\mathrm{d}x)\bigg) 
\end{aligned}\\ \notag
&= \gamma \, \BE[V_i^{(k)}(0,p) V_j^{(l)}(0,p)] + 
\gamma^2 \int \BE[V_i^{(k)}(x,0,p) \, V_j^{(l)}(0,x,p)]
- \BE V_i^{(k)}(0,p) \, \BE V_j^{(l)}(0,p) \mathrm{d}x
\end{align}
and that this asymptotic covariance is finite. In the second step, 
we will show that the asymptotic covariance considered in \eqref{firststep} 
equals $\rho_{i,j}^{k,l}(p)$ (up to a constant), i.e.
\begin{align} \label{asymptCov}
&(d-k+1)(d-l+1) \rho_{i,j}^{k,l}(p) \notag\\
&=
\begin{aligned}[t]
\lim_{t \to \infty} \frac{1}{t} \CV\bigg(&\int V_i^{(k)}(x,p) \, 
\I\{x \in W_t\} \, \eta(\mathrm{d}x),\int V_j^{(l)}(x,p) \, 
\I\{x \in W_t\} \, \eta(\mathrm{d}x)\bigg).
\end{aligned}
\end{align}

We use the Mecke formula, see e.g. \cite{SW}, to get
\begin{align*}
&\CV\bigg(\int V_i^{(k)}(x,p) \, \I\{x \in W_t\} \, 
\eta(\mathrm{d}x),\int V_j^{(l)}(x,p) \, \mathbf{1}\{x \in W_t\} \, \eta(\mathrm{d}x)\bigg)\\
&= \BE\int V_i^{(k)}(x,p) \, V_j^{(l)}(x,p) \, \I\{x \in W_t\} \, \eta(\mathrm{d}x)\\
&\quad+ \BE\iint V_i^{(k)}(x,p) \, V_j^{(l)}(y,p) \, \mathbf{1}\{x \neq y\} \, 
\I\{x,y \in W_t\} \, \eta(\mathrm{d}x) \, \eta(\mathrm{d}y)\\
&\quad- \BE\int V_i^{(k)}(x,p) \, \I\{x \in W_t\} \, \eta(\mathrm{d}x) \, 
\BE\int V_j^{(l)}(x,p)\, \I\{x \in W_t\} \, \eta(\mathrm{d}x)\\
&= \gamma \int \BE[V_i^{(k)}(x,p) \, V_j^{(l)}(x,p)] \, \mathbf{1}\{x \in W_t\} 
\, \mathrm{d}x\\
&\quad+
\begin{aligned}[t]
\gamma^2 \iint \big(&\BE[V_i^{(k)}(x,y,p) \, V_j^{(l)}(y,x,p)]- \BE V_i^{(k)}(x,p) \, 
\BE V_j^{(l)}(y,p) \big) \mathbf{1}\{x,y \in W_t\} \, \mathrm{d}x \, \mathrm{d}y.
\end{aligned}
\end{align*}
Using stationarity of $\eta$, translation invariance of the functions 
$V_i^{(k)}(\cdot,p)$, $V_i^{(k)}(\cdot,\cdot,p)$ and a change of variables we get
\begin{align*}
&\frac{1}{V_d(W_t)}\CV\bigg(\int V_i^{(k)}(x,p) \, \mathbf{1}\{x \in W_t\} \, 
\eta(\mathrm{d}x),\int V_j^{(l)}(x,p) \, \mathbf{1}\{x \in W_t\} \, \eta(\mathrm{d}x)\bigg)\\
&= \gamma \, \BE[V_i^{(k)}(0,p) \, V_j^{(l)}(0,p)]\\
&\quad+ \frac{\gamma^2}{V_d(W_t)} \iint \bigg(
\begin{aligned}[t]
&\vphantom{\int}\BE[V_i^{(k)}(x,0,p) \, V_j^{(l)}(0,x,p)]- \BE V_i^{(k)}(0,p) \, 
\BE V_j^{(l)}(0,p)\bigg)\\
&\times \mathbf{1}\{x+y,y \in W_t\} \, \mathrm{d}x \, \mathrm{d}y
\end{aligned}\\
&= \gamma \, \BE[V_i^{(k)}(0,p) \, V_j^{(l)}(0,p)]\\
&\quad+
\begin{aligned}[t]
\gamma^2 \int \big(&\BE[V_i^{(k)}(x,0,p) \, V_j^{(l)}(0,x,p)]- \BE V_i^{(k)}(0,p) \, 
\BE V_j^{(l)}(0,p)\big) \, \frac{V_d(W_t \cap (W_t-x))}{V_d(W_t)} \, \mathrm{d}x.
\end{aligned}
\end{align*}
By the dominated convergence theorem this converges to the right-hand side of 
\eqref{firststep} as $t \to \infty$, which is in fact the right-hand side of 
\eqref{rhoijkl}. Since $V_d(W_t \cap (W_t-x))/t \leq 1$, a dominating 
function can be given by
\begin{align}\label{Majorante}
|\BE [V_i^{(k)}(x,0,p) \, V_j^{(l)}(0,x,p)]- \BE V_i^{(k)}(0,p) \, \BE V_j^{(l)}(0,p)|.
\end{align}
In the following we will show that the integral of this function is finite.

Next, we will need a technical tool. The \emph{Voronoi flower} of 
$x \in \eta$ is defined by
\begin{align*}
S(\eta,x) := \bigcup_{y \in C(\eta,x)} B(y,\|y-x\|),
\end{align*}
where $B(x,r)$ denotes the closed ball with center $x \in \R^d$ and radius $r\geq 0$. 
Using this definition and the triangle inequality we obtain
\begin{align*}
&\int |\BE[V_i^{(k)}(x,0,p) \, V_j^{(l)}(0,x,p)] - \BE V_i^{(k)}(0,p) \, 
\BE V_j^{(l)}(0,p)| \, \mathrm{d}x\\
&\leq \int
\begin{aligned}[t]
& \bigg| \BE\bigg[
\begin{aligned}[t]
&\I\bigg\{S(\eta^{0,x},x) \subset B\bigg(x,\frac{\|x\|}{3}\bigg)\bigg\} V_i^{(k)}(x,0,p)\\
&\times \I\bigg\{S(\eta^{0,x},0) \subset B\bigg(0,\frac{\|x\|}{3}\bigg)\bigg\} \, 
V_j^{(l)}(0,x,p)\bigg]
\end{aligned}\\
&- \BE V_i^{(k)}(0,p) \, \BE V_j^{(l)}(0,p) \bigg| \, \mathrm{d}x
\end{aligned}\\
&\quad+ \int \bigg| \BE\bigg[
\begin{aligned}[t]
&\bigg(
\begin{aligned}[t]
&\I\bigg\{S(\eta^{0,x},0) \not \subset B\bigg(0,\frac{\|x\|}{3}\bigg)\bigg\} 
+ \I\bigg\{S(\eta^{0,x},x) \not \subset B\bigg(x,\frac{\|x\|}{3}\bigg)\bigg\}\\
&- \I\bigg\{S(\eta^{0,x},0) \not \subset B\bigg(0,\frac{\|x\|}{3}\bigg), 
S(\eta^{0,x},x) \not \subset B(x,\frac{\|x\|}{3})\bigg\}\bigg)
\end{aligned}\\
&\times V_i^{(k)}(x,0,p) \, V_j^{(l)}(0,x,p) \bigg] \bigg| \, \mathrm{d}x
\end{aligned}\\
&=: I_1 + I_2.
\end{align*}
In the following we assume $x \neq 0$. Later we will need that
\begin{align} \label{VF}
\I\bigg\{S(\eta^{0,x},x) \subset B\bigg(x,\frac{\|x\|}{3}\bigg)\bigg\} 
= \I\bigg\{S(\eta^{x},x) \subset B\bigg(x,\frac{\|x\|}{3}\bigg)\bigg\}.
\end{align}
Indeed, in the case $S(\eta^{0,x},x) \subset B(x,\|x\|/3)$ 
the origin cannot be contained in $S(\eta^{0,x},x)$ and is therefore not a 
neighbour of $x$ with respect to $\eta^{0,x}$. Because $S(\eta^{0,x},x)$ is 
determined by $x$ and the neighbours of $x$ with respect to $\eta^{0,x}$, 
the deletion of the origin does not change the Voronoi flower of $x$, 
i.e.\ $S(\eta^x,x)=S(\eta^{0,x},x) \subset B(x,\|x\|/3)$. 
In the case $S(\eta^{0,x},x) \not\subset B(x,\|x\|/3)$, the Voronoi 
flower cannot get larger if we add more points to the point process, 
i.e.\ $S(\eta^{0,x},x) \subset S(\eta^x,x)$. This implies 
$S(\eta^x,x)\not\subset B(x,\|x\|/3)$.

Now we will use the stopping set property of the Voronoi 
flowers $S(\eta^x,x)$ and $S(\eta^0,0)$, see \cite{Zu99}. 
Because the Voronoi cell and the corresponding Voronoi flower are 
determined by the Poisson points contained in the flower, the random variable
$$ 
\I\bigg\{S(\eta^{x},x) \subset B\bigg(x,\frac{\|x\|}{3}\bigg)\bigg\} \, V_i^{(k)}(x,0,p) 
$$
is determined by the intersection of $\eta$ and $B(x, \|x\|/3)$. Analogously, 
the random variable
$$ 
\I\bigg\{S(\eta^{0},0) \subset B\bigg(0,\frac{\|x\|}{3}\bigg)\bigg\} \, V_j^{(l)}(0,x,p) 
$$
is determined by the intersection of $\eta$ and $B(0, \|x\|/3)$. Using \eqref{VF}, 
an analogous equation for the Voronoi flower of the origin and the 
fact that $B(x,\|x\|/3)$ and $B(0, \|x\|/3)$ are disjoint and the restrictions of 
a Poisson process to disjoint sets are independent, we have
\begin{align*}
I_1 &=
\begin{aligned}[t]
\int \bigg| & \BE\bigg[\I\bigg\{S(\eta^{x},x) \subset B\bigg(x,\frac{\|x\|}{3}\bigg)\bigg\} 
\, V_i^{(k)}(x,p)\bigg]\\
&\times \BE\bigg[\I\bigg\{S(\eta^{0},0) 
\subset B\bigg(0,\frac{\|x\|}{3}\bigg)\bigg\} \, V_j^{(l)}(0,p)\bigg] - \BE V_i^{(k)}(0,p) \, \BE V_j^{(l)}(0,p) \bigg| \, \mathrm{d}x.
\end{aligned}
\end{align*}
Thus, we get by the stationarity
\begin{align*}
I_1 &= \int
\begin{aligned}[t]
\bigg| & \BE\bigg[\I\bigg\{S(\eta^{0},0) \subset B\bigg(0,\frac{\|x\|}{3}\bigg)\bigg\} 
\, V_i^{(k)}(0,p)\bigg]\\
&\times \BE\bigg[\I\bigg\{S(\eta^{0},0) \subset 
B\bigg(0,\frac{\|x\|}{3}\bigg)\bigg\} \, V_j^{(l)}(0,p)\bigg]- \BE V_i^{(k)}(0,p) \, \BE V_j^{(l)}(0,p) \bigg| \, \mathrm{d}x
\end{aligned}\\
&= \int
\begin{aligned}[t]
\bigg| & \BE\bigg[\I\bigg\{S(\eta^0,0) \not \subset B\bigg(0,\frac{\|x\|}{3}\bigg)\bigg\} 
\, V_i^{(k)}(0,p)\bigg]\\
&\times \BE\bigg[\I\bigg\{S(\eta^0,0) \not \subset B\bigg(0,\frac{\|x\|}{3}\bigg)\bigg\} \, V_j^{(l)}(0,p)\bigg]\\
&-\BE\bigg[\I\{S(\eta^0,0) \not \subset B\bigg(0,\frac{\|x\|}{3}\bigg)\bigg\} \, 
V_i^{(k)}(0,p)\bigg] \, \BE V_j^{(l)}(0,p)\\
&- \BE V_i^{(k)}(0,p) \, \BE\bigg[\I\bigg\{S(\eta^0,0) \not 
\subset B\bigg(0,\frac{\|x\|}{3}\bigg)\bigg\} \, V_j^{(l)}(0,p)\bigg] \bigg| \, \mathrm{d}x.
\end{aligned}
\end{align*}
After using the triangle inequality the first summand is smaller than the 
second and the second and the third summand are the same (up to different parameters), 
so it's enough to show the finiteness of the second. Using the definition of 
neighbourhood we have
\begin{align*}
&\int \BE\bigg[\I\bigg\{S(\eta^0,0) \not \subset B\bigg(0,\frac{\|x\|}{3}\bigg)\bigg\} 
\, V_i^{(k)}(0,p)\bigg] \, \BE V_j^{(l)}(0,p) \, \mathrm{d}x\\
&\leq \int
\begin{aligned}[t]
&\BE\bigg[\I\bigg\{2 \, \mathrm{diam}(C(\eta^0,0)) > \frac{\|x\|}{3}\bigg\} 
\, |N(\eta^0,0)|^{d-k} \, V_i(C(\eta^0,0))\bigg]\\
&\times \BE[|N(\eta^0,0)|^{d-l} \, V_j(C(\eta^0,0))] \, \mathrm{d}x
\end{aligned}\\
&\leq
\begin{aligned}[t]
&\int \BP\bigg(\mathrm{diam}(C(\eta^0,0)) > \frac{\|x\|}{6}\bigg)^{1/3} \, \mathrm{d}x\\
&\times (\BE |N(\eta^0,0)|^{3d-3k} \, \BE V_i(C(\eta^0,0))^3)^{1/3} \, 
(\BE |N(\eta^0,0)|^{2d-2l} \, \BE V_j(C(\eta^0,0))^2)^{1/2}.
\end{aligned}
\end{align*}
This is finite because of Lemma \ref{Nachbarzahl} and Corollary \ref{Momente}. 
Because $i,j,k,l$ were arbitrary we have shown the finiteness of the third summand as well.

To show the finiteness of $I_2$ it is again enough to show the finiteness of
 the first summand. We have
\begin{align*}
&\int \BE\bigg[\I\bigg\{S(\eta^{0,x},0) \not \subset B\bigg(0,\frac{\|x\|}{3}\bigg)\bigg\} 
\, V_i^{(k)}(x,0,p) \, V_j^{(l)}(0,x,p) \bigg] \, \mathrm{d}x\\
&\leq \int \BE\bigg[
\begin{aligned}[t]
&\I\bigg\{S(\eta^{0,x},0) \not \subset B\bigg(0,\frac{\|x\|}{3}\bigg)\bigg\}\\
&\times |N(\eta^{0,x},x)|^{d-k} \, V_i(C(\eta^{0,x},x)) \, |N(\eta^{0,x},0)|^{d-l} 
\, V_j(C(\eta^{0,x},0))\bigg] \, \mathrm{d}x
\end{aligned}\\
&\leq \int \BP\bigg(\mathrm{diam}(C(\eta^0,0))> \frac{\|x\|}{6}\bigg)^{1/5} \, \mathrm{d}x\\
&\quad \times (\BE (|N(\eta^0,0)|+1)^{5d-5k} \, \BE (|N(\eta^0,0)|+1)^{5d-5l} \, 
\BE V_i(C(\eta^0,0))^5 \, \BE V_j(C(\eta^0,0))^5)^{1/5}.
\end{align*}
This is finite because of Lemma \ref{Nachbarzahl}, 
Lemma \ref{expAbfdiam} and Corollary \ref{Momente}. Hence, the integral 
of the dominating function given in \eqref{Majorante} is finite and therefore 
the second summand of the right-hand side of \eqref{firststep} is finite, too.

Using the monotonicity of the intrinsic volumes, normality and H\"older's 
inequality we get for the first summand of the right-hand side of \eqref{firststep}
\begin{align*}
\BE[V_i^{(k)}(0,p) \, V_j^{(l)}(0,p)] &\leq \BE[|\cF_k(C(\eta^0,0))| \, V_i(C(\eta^0,0)) 
\, |\cF_l(C(\eta^0,0))| \, V_j(C(\eta^0,0))]\\
&\leq \BE[|N(\eta^0,0)|^{d-k} \, V_i(C(\eta^0,0)) \, |N(\eta^0,0)|^{d-l} \, V_j(C(\eta^0,0))]\\
&\leq (\BE |N(\eta^0,0)|^{6d-3k-3l} \, \BE V_i(C(\eta^0,0))^3 \, \BE V_j(C(\eta^0,0))^3)^{1/3}.
\end{align*}
This is finite by Lemma \ref{Nachbarzahl} and Corollary \ref{Momente}. So, 
the right-hand side of \eqref{firststep} is finite.

The next step is to prove \eqref{asymptCov}. Because of the normality of the 
Poisson Voronoi tessellation we have
\begin{align}\label{umschreiben}
&\int V_i(F(x) \cap W_t) \, f_n^k(|\cS_n(x)|,p) \, \eta^{(k)}(\mathrm{d}x) \notag\\
&= \frac{1}{d-k+1} \int V_i(\cF_k(C(\eta,x)) \cap W_t) \, f_n^k(|\cS_n(x)|,p) \, 
\eta(\mathrm{d}x)
\end{align}
and \eqref{asymptCov} is equivalent to
\begin{align}\label{equiv}
&\begin{aligned}[t]
\lim_{t \rightarrow \infty} \frac{1}{t} \CV\bigg(&\int V_i(\cF_k(C(\eta,x)) \cap W_t) 
\, f_n^k(|\cS_n(x)|,p) \, \eta(\mathrm{d}x),\\
&\int V_j(\cF_l(C(\eta,x)) \cap W_t) f_n^l(|\cS_n(x)|,p) \, \eta(\mathrm{d}x) \bigg)
\end{aligned} \notag\\
&= \lim_{t \rightarrow \infty} \frac{1}{t} \CV\bigg( \int V_i^{(k)}(x,p) 
\, \I\{x \in W_t\} \, \eta(\mathrm{d}x), \int V_j^{(l)}(x,p) \, 
\I\{x \in W_t\} \, \eta(\mathrm{d}x) \bigg).
\end{align}
Analogously to \eqref{7.3a} we obtain that
\begin{align}\label{7.3}
\lim_{t \to \infty} \frac{1}{t} \BV\bigg(\int V_i(\cF_k(C(\eta,x)) \cap W_t) \, 
f_n^k(|\cS_n(x)|,p) - V_i^{(k)}(x,p) \, \mathbf{1}\{x \in W_t\} \, 
\eta(\mathrm{d}x) \bigg) =0.
\end{align}

By using \eqref{7.3} and the Cauchy-Schwarz inequality we will show that \eqref{equiv} 
(and therewith \eqref{asymptCov}) holds. We abbreviate
\begin{align*}
V_i^{k} &:= \int V_i(\cF_k(C(\eta,x)) \cap W_t) \, f_n^k(|\cS_n(x)|,p) \, \eta(\mathrm{d}x),\\
U_i^{k} &:= \int V_i(\cF_k(C(\eta,x))) \, \mathbf{1}\{x \in W_t\} \, f_n^k(|\cS_n(x)|,p) \, 
\eta(\mathrm{d}x).
\end{align*}
Note that
\begin{align*}
&\CV(V_i^{k},V_j^{l})-\CV(U_i^{k},U_j^{l})\\
&= \CV(V_i^{k}-U_i^{k},U_j^{l}) + \CV(V_i^{k}-U_i^{k},V_j^{l}-U_j^{l}) 
+ \CV(U_i^{k},V_j^{l}-U_j^{l}).
\end{align*}
Hence we obtain by the Cauchy-Schwarz inequality
\begin{align*}
0 &\leq \lim_{t \to \infty} \frac{1}{t} |\CV(V_i^{k},V_j^{l})-\CV(U_i^{k},U_j^{l})|\\
&\leq
\begin{aligned}[t]
\lim_{t \to \infty} &\frac{1}{t} |\CV(V_i^{k}-U_i^{k},U_j^{l})| 
+ \frac{1}{t} |\CV(V_i^{k}-U_i^{k},V_j^{l}-U_j^{l})| + \frac{1}{t} |\CV(U_i^{k},V_j^{l}-U_j^{l})|
\end{aligned}\\
&\leq
\begin{aligned}[t] 
\lim_{t \to \infty} &\sqrt{\frac{1}{t} \BV(V_i^{k}-U_i^{k})} \sqrt{\frac{1}{t} \BV(U_j^{l})}
 + \sqrt{\frac{1}{t} \BV(V_i^{k}-U_i^{k})} \sqrt{\frac{1}{t} \BV(V_j^{l}-U_j^{l})}\\
&+ \sqrt{\frac{1}{t} \BV(U_i^{k})} \sqrt{\frac{1}{t} \BV(V_j^{l}-U_j^{l})}.
\end{aligned}
\end{align*}
This is zero because of \eqref{7.3} and the already shown finiteness of 
$\lim_{t \to \infty} \frac{1}{t} \BV(U_i^{k})$.
\qed

\bigskip
The variance of the Euler characteristic of a planar
Poisson Voronoi percolation is worth special mentioning:

\begin{corollary}\label{varianceEuler} Consider cell percolation
on a planar Poisson Voronoi tessellation. Then
the asymptotic variance $\sigma_{0,0}$ of the Euler characteristic
exists and is given by
\begin{align*}
\sigma_{0,0}(p)= \gamma_2 \mu_2 p^3(1-p)^3 
+ \gamma_2 p(1-p)(1-8p-6p^2+28p^3-14p^4).
\end{align*}
Moreover, $\sigma_{0,0}$ has a strict global maximum at $1/2$.
\end{corollary}
{\em Proof.} By Theorem \ref{thm7.1} we are allowed to apply
Theorem \ref{Kovarianzstruktur}. Since the asymptotic variance $\tau_{0,0}^{2,2}$ 
equals the intensity $\gamma_2$ we obtain the formula for
$\sigma_{0,0}$. The second assertion follows from the corresponding
assertion of Theorem \ref{Kovarianzstruktur} and \eqref{36}
or from a direct calculation.\qed

\begin{appendix}
\section{Appendix}

We consider a Poisson Voronoi tessellation $X$ 
generated by Poisson process  $\eta$ of intensity $\gamma>0$. 
Recall the definition \eqref{Nachbarschaft2} 
of the neighbourhood of second order.

\begin{lemma}\label{Nachbarzahl}
All moments of the cardinality of the neighbourhood of second order of the 
typical cell are finite, that is
$\BE |N_2(\eta^0,0)|^m < \infty$ for all $m \in \N$.
\end{lemma}

{\em Proof.}
Let $X_1,X_2,\ldots$ be the enumeration of $\eta$ such that $0<\|X_1\|<\|X_2\|<\ldots$. 
Let $n:=|N_2(\eta^0,0)|$ and assume $n\geq 2$ because this is almost surely 
satisfied. Hence, there is a point $x \in N_2(\eta^0,0)$ with $\|x\| \geq \|X_n\| =:t$. 
It is easy to see that $x \in \eta$ is a neighbour of second order of $0$ with 
respect to $\eta^0$ if and only if there exist $y \in \eta \setminus \{x\}$ and 
balls $B$ and $B'$ with $0,y \in \partial B$, $x,y \in \partial B'$, 
$\mathrm{int}(B) \cap \eta =\emptyset$ and $\mathrm{int}(B') \cap \eta =\emptyset$. 
Since $\|x\| \geq t$ we have $\mathrm{diam}(\tilde{B}) \geq t/2$ and 
$\tilde{B} \cap B(0,t/2) \neq \emptyset$ for either $\tilde{B}=B$ or $\tilde{B}=B'$.

There exist balls $B_1,\ldots,B_l \subset \mathrm{int}(B(0,t))$ with diameter $t/8$ 
such that each ball $\tilde{B}$ with diameter at least $t/2$ and 
$\tilde{B} \cap B(0,t/2) \neq \emptyset$ contains at least one of the balls 
$B_1,\ldots,B_l$. By a scaling argument, the number $l$ of balls can be 
chosen independently of $t$.

So, $|N_2(\eta^0,0)|\geq n$ implies $\eta(B_i)=0$ for at least one $i \in \{1,\ldots,l\}$. 
Because of the binomial property of the Poisson process we have
\begin{align*}
\BP(\eta(B_1)=0|\|X_n\|=t) &= \bigg(1-\frac{\kappa_d(t/16)^d}{\kappa_dt^d}\bigg)^{n-1}
= \bigg(1-\frac{1}{16^d}\bigg)^{n-1}
\end{align*}
and the conditional probability that at least one ball $B_i$ contains no point of $\eta$ 
is at most $l(1-1/16^d)^{n-1}$. If we denote the density of $\|X_n\|$ by $f_n$, we get
\begin{align*}
\BP(|N_2(\eta^0,0)| \geq n) = \int \BP(|N_2(\eta^0,0)| \geq n| \|X_n\|=t) \, f_n(t) \, 
\mathrm{d}t \leq \int l\bigg(1-\frac{1}{16^d}\bigg)^{n-1} f_n(t) \, \mathrm{d}t.
\end{align*}
Because $\kappa_d \|X_n\|^d$ is gamma distributed with parameters $n$ and 1, 
see also \cite{BaumLa07}, we get from an obvious transformation
\begin{align*}
f_n(t) = \frac{\kappa_d^n d}{(n-1)!} t^{dn-1} \exp(-\kappa_dt^d).
\end{align*}
A change of variables gives
\begin{align}\label{Nach2tail}
\BP(|N_2(\eta^0,0)| \geq n) &\leq \frac{l(1-1/16^d)^{n-1}\kappa_d^nd}{(n-1)!} 
\int t^{dn-1} \exp(-\kappa_dt^d) \, \mathrm{d}t \notag\\
&= \frac{l(1-1/16^d)^{n-1}\kappa_d^n}{(n-1)!} \int s^{n-1} \exp(-\kappa_ds) \, 
\mathrm{d}s = l(1-1/16^d)^{n-1}.
\end{align}
This implies the assertion.
\qed \bigskip

\begin{lemma}\label{Appendix2}
All moments of the cardinality of neighbouring $k$-faces of the typical
$n$-face in a 
Poisson Voronoi tessellation exist, i.e.
\begin{align*}
\sum_{k=0}^d \BE_n^0 |\cS_k(0)|^m < \infty, \quad m \in \N.
\end{align*}
\end{lemma}

\emph{Proof.} 
By normality it suffices to treat the case
$n>k$. Proposition \ref{t2.1} implies
\begin{align*}
(d-n+1) \, \gamma_n \BE_n^0 |\cS_k(0)|^m 
&= \gamma_n \, \BE_n^0\bigg[\sum_{F \in \cS_d(0)} |\cS_k(0)|^m \bigg] 
= \gamma_d \, \BE_d^0\bigg[\sum_{F \in \cS_n(0)} |\cS_k(s(F))|^m\bigg].
\end{align*}
Because of $n>k$ we can bound $|\cS_k(s(F))|$ for $F \in \cS_n(0)$ by $|\cS_k(0)|$. 
Further, each $k$-face of the typical cell is contained in exactly $d-k$ 
neighbours of the typical cell and we get
\begin{align*}
(d-n+1) \, \gamma_n \BE_n^0 |\cS_k(0)|^m 
&\leq \gamma_d \BE_d^0[|\cS_n(0)| \, |\cS_k(0)|^m] 
\leq  \gamma \, \BE^0[|N(\eta^0,0)|^{d-n+md-mk}]
\end{align*}
and this is finite because of Lemma \ref{Nachbarzahl}.
\qed


\bigskip

The proof of the next lemma is given in \cite{HugSchn07}, Theorem 2.

\begin{lemma}\label{expAbfdiam}
There are constants $c_1,c_2>0$ such that for all $u \in \R$
\begin{align*}
\BP(\mathrm{diam}(C(\eta^0,0)) \geq u) \leq c_1 \exp(-c_2u).
\end{align*}
\end{lemma}

\begin{corollary}\label{Momente}
All moments of the intrinsic volumes of the typical cell are finite, i.e. for $i \in \{0,\ldots,d\}$
\begin{align*}
\BE V_i(C(\eta^0,0))^m < \infty, \quad m \in \N.
\end{align*}
\end{corollary}

{\em Proof.}
We use $C(\eta^0,0) \subset B(0,\mathrm{diam}(C(\eta^0,0)))$ and the monotonicity of the intrinsic volumes to get for any $m \in \N$
\begin{align*}
\BE V_i(C(\eta^0,0))^m &\leq \BE V_i(B(0,\mathrm{diam}(C(\eta^0,0))))^m = V_i(B(0,1))^m \, \BE \, \mathrm{diam}(C(\eta^0,0))^{im}.
\end{align*}
This is finite because of Lemma \ref{expAbfdiam}.
\qed \bigskip

We will introduce a modification of the system $\cS_l(\varphi,x)$ for a tessellation $\varphi \in \mathbf{T}$ and $x \in \R^d$ with $F(x) \in \cF_k(\varphi)$ and define
\begin{align}\label{VerallgSeitenstern}
\tilde{\cS}_l(\varphi,x) :=
\begin{cases}
\cS_l(\varphi,x), \quad \min(k,l)<d,\\
\{G \in \cF_k(\varphi): G \cap F(x) \in \cF_{k-1}(\varphi)\}, \quad k=l=d.
\end{cases}
\end{align}
We will use the following exchange formula, which is similar to Proposition \ref{t2.1}, but considers the system $\tilde{\cS}_l(0)$ instead of $\cS_l(0)$. The proof can be easily given by Neveu's exchange formula.

\begin{proposition}\label{t2.1a}
Let $k,l \in \{0,\ldots,d\}$ and $g : \cP^d \times \cP^d \rightarrow [0,\infty)$ be a measurable function. Then
\begin{align}\label{2.18a}
&\gamma_k \BE_k^0 \sum_{G \in \tilde{\cS}_l(0)} g(F(0), G-s(G)) = \gamma_l \BE_l^0 \sum_{F \in \tilde{\cS}_k(0)} g(F-s(F), F(0)).
\end{align}
\end{proposition}

Let $R(B)$ be the radius of the circumball of a subset $B \subset \R^d$.

\begin{lemma}\label{Umkugelradius2}
We have for all $m \in \N$
\begin{align*}
\sum_{k=0}^d \BE_k^0 \bigg[R\bigg(F(0) \cup \bigcup_{G \in \tilde{\cS}_d(0)}G\bigg)^m \bigg] < \infty.
\end{align*}
\end{lemma}

{\em Proof.}
For $k<d$ we have
\begin{align*}
\BE_k^0 R\bigg(F(0) \cup \bigcup_{G \in \tilde{\cS}_d(0)}G\bigg)^m \leq \BE_k^0 \max_{G \in \tilde{\cS}_d(0)} (2R(G))^m \leq 2^m \BE_k^0 \sum_{G \in \tilde{\cS}_d(0)} R(G)^m
\end{align*}
and in the case $k=d$
\begin{align*}
\BE_d^0 R\bigg(F(0) \cup \bigcup_{G \in \tilde{\cS}_d(0)}G\bigg)^m &\leq \BE_d^0 (R(F(0)) +2 \max_{G \in \tilde{\cS}_d(0)} R(G))^m\\
&\leq 2^m \, \BE_d^0 \bigg[R(F(0))^m+ 2^m \sum_{G \in \tilde{\cS}_d(0)} R(G)^m\bigg].
\end{align*}
Because of $R(F(0)) \leq 2\, \mathrm{diam}(C(\eta^0,0))$ and Lemma \ref{expAbfdiam} it is in both cases enough to show that
$$ \BE_k^0 \sum_{G \in \tilde{\cS}_d(0)} R(G)^m < \infty. $$

If we define $g:\cP^d \times \cP^d \rightarrow [0,\infty)$ by $g(F,G):= R(G)^m$, \eqref{2.18a} and the Cauchy-Schwarz inequality yield
\begin{align*}
\gamma_k \BE_k^0 \sum_{G \in \tilde{\cS}_d(0)}R(G)^m = \gamma_d \BE_d^0 [|\tilde{\cS}_k(0)| \, R(F(0))^m] \leq \gamma_d \BE_d^0[|\tilde{\cS}_k(0)|^2]^{1/2} \, \BE_d^0[R(F(0))^{2m}]^{1/2}.
\end{align*}
As above, the second factor is finite. For $k<d$ it follows by normality that $|\tilde{\cS}_k(0)| = d-k+1$ and for $k=d$ we have $|\tilde{\cS}_k(0)| \leq |N_2(\eta^0,0)|$, so the first factor is finite because of Lemma \ref{Nachbarzahl}.
\qed

\end{appendix}

\end{document}